\documentclass[11pt,a4paper]{article}
\usepackage{amsmath}
\usepackage{amssymb,amsfonts, amstext,amsmath}
\usepackage[latin1]{inputenc}
\usepackage{textcomp}
\usepackage{theorem}

\usepackage{hyperref}

\usepackage{algorithm2e}
\SetKw{KwStep}{step}

\usepackage{graphicx}
\usepackage{xcolor}
\usepackage{fancybox}

\newcounter{exo}

\newcommand{\Q}{\mathbb{Q}}

\newcommand{\R}{\mathbb{R}}
\newcommand{\Z}{\mathbb{Z}}
\newcommand{\N}{\mathbb{N}}

 \newcommand{\tend}{\xrightarrow[n \to \infty]{}}

\newcommand{\CFE}{\mathcal{C}}

\renewcommand{\phi}{\varphi}

 \newcommand{\esp}{\hspace*{0.3cm}}

\newcommand{\titre}[1]{\begin{center}  \textbf{\Large{#1}} \\ \end{center}}

\newcommand{\resultat}[1]{  \shadowbox{\parbox{16cm}{#1}}}

\newcommand{\inter}[2]{ \{  #1 \cdots #2 \}  }

 \theorembodyfont{\sl}  \newtheorem{montheo}{Theorem}
 \theorembodyfont{\normalfont}  
     \newtheorem{madef}{Definition} 
       \newtheorem{monalgo}{Algorithm} 
    
      \newtheorem{maprop}{Proposition} 
      \newtheorem{moncoro}{Corollary}     
                              
 \theorembodyfont{\sl} \newenvironment{mademo}{\textbf{Proof : } \small \\}{\normalsize $\blacksquare$ \\ }
 \newtheorem{monlem}{Lemma}

\begin {document}
 
 \titre{A variant of Ostrowski numeration }
 
\hspace{6cm} \large Emmanuel Cabanillas  \normalsize \\

\vspace*{1cm}

\hspace{7cm} ABSTRACT : \\
 
\esp In this article, we propose a variant of the usual Ostrowski $\alpha$-numeration ( where $\alpha$ is a real in $[0,1[$) that codes  integers ( positive as well as negative) and reals of $[0,1[$ ( instead of $[-\alpha, 1-\alpha[$), so that for every  integer $n$, $n$ and $\{ n\alpha \}$ have the same coding sequence. These coding sequences respect natural lexicographic orders and will be used to prove well known results on order properties of Kronecker sequences $(\{ n\alpha - \beta \})_n$.  \\

\esp

\vspace*{1.5cm}

 \tableofcontents

 \newpage
 
 \section{Introduction}
 
   \subsection{overview}
   
   \esp Ostrowski's numeration system is based on convergents $(q_n)_{n\in\N}$ of a real $\alpha\in[ 0,1[$ and code, with a sequence of digits non negative integers as well as reals in $[-\alpha,1-\alpha[ $ ( see \textbf{[6]}  for the original article and \textbf{[1]} for a survey). Definitions are mentioned in \ref{sec:ost}\\
   \esp In \ref{sec:rational-numeration} and \ref{sec:irrational-numeration}, we propose a variant of this system : it is still based on  $(q_n)_n$, but the " markovian condition" is changed and we will be able to code any integer $n$ and any real $\{ n\alpha \}$ with the same finite sequence ( $\{ x \}$ denotes the fractional part of a real $x$). We study separately the cases $\alpha$ irrational and $\alpha$ rational. This last case could appear uninteresting, but it is useful for applications to numerical semigroups for example ( see \textbf{[3]}). \\
     \esp In \ref{ch:dyn}, we give some dynamical aspects of this $\alpha$-numeration. \\
     \esp In \ref{ch:kro}, we use it to explore some order properties of   Kronecker sequences $(\{ n\alpha +\beta\})_n$, as the famous " three distance theorem". These sequences have been widely studied with various points of view and we refer to \textbf{[1]} for an exhaustive bibliography.

   \subsection{notations}

   \esp All along this paper, we will denote :  $\Z$ the set of integers, $\N^*$ the set of positive integers and  $\N$ the set of non negative integers.\\
   \esp For all reals $x$, $\lfloor x \rfloor$ denotes its floor ,$\lceil x \rceil $ its ceiling and $\{ x \}$ its fractional part.\\
   \esp For a sequence $d=(d_k)_{k\in \N^*}$, we use the following notations for slices of $d$ : for all integers $r,s$ such that $0<r \leqslant s$ :
   $$ d_{[r,s]}=(d_r,d_{r+1},\cdots,d_s) \esp ; \esp d_{[r,\infty]}=(d_r,d_{r+1},\cdots)$$
   \esp We will also use concatenation of sequences and intuitive notations as $(3,5,0^4,1,6,0^{\infty})$ to denote $(3,5,0,0,0,0,1,6,0,0,0,\cdots)$. Moreover, if $(a_k)_{k\in\N^*}$ is a sequence of positive integers and if we restrict ourself to sequences in $\prod_k\inter{0}{a_k}$, then $\max$ at the index $k$ will denote $a_k$ : for example, $(\max,1,0,\max,3,\cdots,)$ means $(a_1,1,0,a_4,3,...)$. So, the notation $\max^r$ or $(\max,0)^r$, where $r\in \N\cup \{ \infty \}$ will often be used. For example : $(0^2,\max^3,0^4,(\max,0)^{\infty})$ denotes the sequence $(0,0,a_3,a_4,a_5,0,0,0,0,a_{10},0,a_{12},0,a_{14},0,\cdots)$.\\
    
   \esp For $\alpha$-numeration, we will often use two lexicographic orders on sequences  of $\R^{\N^*}$ : \\
   $\blacktriangleright$ the \textbf{reversed lexicographic order ( RLO) } denoted $ \underset{R}{\leqslant}$  : 
     $$ d \underset{R}{\leqslant} d' \Leftrightarrow d=d' \text{ or } \exists j\in \N^*, \begin{cases} d_j<d'_j \\
     \forall i>j, d_i=d'_i \end{cases} $$
    $\blacktriangleright$ the \textbf{alternate lexicographic order ( ALO)} denoted $ \underset{A}{\leqslant}$ : 
     $$ d \underset{A}{\leqslant}d' \Leftrightarrow d=d' \text{ or } \exists j\in \N^*, \begin{cases} (-1)^{j-1}d_j<(-1)^{j-1}d'_j \\
     \forall i\in \inter{1}{j-1}, d_i=d'_i \end{cases} $$
  \esp  ALO is a total order on $\R^{\N^*}$, but RLO is only a partial order on $\R^{\N^*}$. Now, RLO is a total order on on $\R^{(\N^*)}$, the set of real sequences that ends with $0^{\infty}$.\\
  \esp We will also use ALO with a shift on indices   for continued fraction expansions in 1.3 ( named \textbf{CFE} in this paper).
  
  \newpage
 
 \subsection{continued fraction expansions}

     \esp All results given in this subsection are well known and we just want to underline some notations and simple facts.     \\

%     $\bullet$ For all $x,y\in \R$, with $y\not = 0$, we denote $[x,y]=x+\frac{1}{y}$ and extend this notation to finite sequences, by " right associativity" : $[x,y,z]=[x,[y,z]]$.
   $\bullet$ Every irrational $\theta$ can be uniquely represented by its continued fraction expansion ( CFE) and we will write $\theta = [t_0,t_1,\cdots ]=[t_k]_{k\in\N}$, such that $t_k\in\N^*$ for all $k\in\N^*$ and $t_0\in\Z$. $\theta$ is the limit of the " convergents" $([t_0,t_1,\cdots, t_n])_n$, a  sequence of rationals defined inductively by :
   $$ \forall x \in \R, \forall x_1,\cdots,x_n \in \R_+^*,[x]=x \esp ; \esp [x,x_1,\cdots,x_n]=x+\frac{1}{[x_1,\cdots,x_n]} 
   \esp (1)$$ 
   \esp We will denote, for all integer $n$, $\frac{p_n}{q_n}$ ( or $\frac{p_n(\theta)}{q_n(\theta)}$ if necessary) the reduced fraction that represents $[t_0,t_1,\cdots, t_n]$.\\
   \esp In addition, if we define $\phi$ : 
   $$ \phi : \begin{cases} \Z\times (\N^*)^{\N^*} \to \R\backslash \Q \\ (t_k)_{k\in\N} \to [t_k]_{k\in\N} \end{cases}$$
   \esp this map is bijective and increasing, with the Alternate Lexicographic Order ( ALO) on $\Z\times (\N^*)^{\N}$ defined by : 
   $$ (t_k)_{k\in\N}\leqslant_A (t'_k)_{k\in\N} \Leftrightarrow  (\forall k\in\N, t_k=t'_k ) \text{ or }  \exists j\in\N, \begin{cases} \forall k\in \inter{0}{j-1}, t_k=t'_k  \\ (-1)^jt_j< (-1)^j t'_j \end{cases}$$
   \esp We also have an expression for the inverse function of $\phi$ : 
   $$  \phi^{-1 } : \begin{cases} \R\backslash \Q \to \Z\times (\N^*)^{\N^*} \\ \theta \to (t_k)_{k\in\N} ,\text{ with } t_0=\lfloor \theta \rfloor ; \forall k\in\N^*, t_k= AT^{k-1}(\{\theta\})\end{cases}$$
   \esp where $T$ is the Gauss map : $]0,1[\to [0,1[, x\to \{ 1/x \}$ and $A : x\to \lfloor 1/x \rfloor$. We know that $T^k(\alpha)\not = 0$ for all $k\in\N$ if and only if $\alpha$ is irrational in $]0,1[$. \\ 
   
   $\bullet$  The case of rationals seems easier, since these one are represented by finite CFE, namely the convergents of irrationals. But, we would like to associate to them infinite CFE, in order to extend $\phi$ to an increasing map with  ALO. \\
   \esp We introduce an $\infty$ number : $\overline{\N^*}$ will denote $\N^*\cup \{ \infty \}$, with the usual extension of the order ( $\forall n\in\N^*, n< \infty$) and of the operations ( $\forall n\in\overline{\N} , n+\infty=\infty$ and $1/\infty=0$). Then, we can end  CFE of rationals with an infinite sequence of $\infty$. With those conventions, the former map $\phi$ extends to an increasing and bijective map $\tilde{\phi}$ from  a subset $E$ of $ \Z\times (\overline{\N^*})^{\N^*}$ to $\R$. Then, $\tilde{\phi}^{-1}$ is given by the same expressions, if we extend $T$ and $A$ to $[0,1[$, with $T(0)=0$ and $A(0)=\infty$.\\
   \esp  We can precise $E$ : it is the set of sequences $(t_k)_k$ such that $t_0\in \Z$ and $t_k\in \overline{\N^*}$ for $k\in\N^*$, such that $t_k=\infty \Rightarrow ( t_{k+1}=\infty$ and ( $t_{k-1}\not = 1$ or $k=1$)). So to say : if the sequence contains $\infty$, the last " finite digit" in the CFE is greater or equal to $2$. We will prefer an alternative way : we will end  CFE of rationals with $[1,\infty^{\infty}]$, where $\infty^{\infty}$ denotes an infinite sequence of $\infty$. Then, we extend naturally the ALO to sequences of CFE, described by :
   
    $$ \CFE= \{ (t_k)\in \Z\times \N^* \times (\overline{\N^*})^{\N}, \forall k\geqslant 2, ( t_k= \infty \Rightarrow ( t_{k+1}=\infty \text{ and } t_{k-1}\in \{ \infty  ,1 \} ) \} $$
   
   \esp The extension of $\phi$ to an increasing and bijective map $\phi_1$ from $\CFE$ to $\R$ is quite natural, but its inverse function will use more complicated maps $T_1$ and $A_1$ : 
 
   \esp We consider the map $I : u\to \lceil u\rceil -1$ and $A_1, T_1$ both defined on $[0,1]$ by :
   $$ A_1 : \begin{cases} 0\to \infty \\ 1\to 1 \\ x\to I(1/x) \text{ if } x \not = 0,1 \end{cases} \esp ; \esp 
    T_1 :\begin{cases} 0\to 0 \\ 1\to 0 \\ x\to 1 \text{ if } 1/x \in \N\backslash \{ 0,1\}\\ x\to \{ 1/x \} \text{ else} \end{cases}$$
    \esp We can now express the inverse function of $\phi_1$ :
      $$   \phi_1^{-1}: \begin{cases} \R\to \CFE \\
             \theta \to (t_k)_{k\in\N} ,\text{ with } t_0=I( \theta ) ; \forall k\in\N^*, t_k= A_1T_1^{k-1}(\theta - I(\theta))\end{cases}$$ 
   \esp For convenience, we abreviate CFE of rationals and omit $\infty^{\infty}$, the infinite " $\infty$" ending sequence. So, $9/4=[2,3,1]$ and $\forall n\in\Z, n=[n-1,1]$.\\

  \textbf{N.B : all along this paper, CFE of a real ( so for any rational) $\alpha$ will denote $\phi_1^{-1}(\alpha)$, but the notation $[t_0,t_1,\cdots,t_k]$ will be more general} ( see (1)).\\

   \subsection{semi-convergents and best rationals}

     $\bullet$ Let $\alpha$ be a real with CFE $[a_k]_{k\in\N}$ and $(p_k/q_k)_k$ its convergents sequence, such that $p_k/q_k=[a_0,\cdots,a_k]$, for all $k$ such that $a_k<\infty$ ( see beginning of this section).\\
     \esp  A \emph{semi-convergent} of $\alpha$ is any rational of the form $\frac{mp_k+p_{k-1}}{mq_k+q_{k-1}}$, with $m\in \inter{0}{a_k}$ and $k\in \N$ such that $a_k<\infty$ ( we take $m>0$ if $k=0$ to avoid $1/0$ !). So, convergents are particular semi-convergents. 
   
   \begin{monlem} Let $\alpha$ be a real with CFE $[a_k]_{k\in\N}$. Semi-convergents of $\alpha$ are exactly the rationals with CFE $[a_0,\cdots,a_{s-1},b_s,1]$, such that $s\in\N,b_s\in\inter{1}{a_s}$ and $a_{s+1}<\infty$.           
    \end{monlem}         
     
     \begin{mademo}
  %  (i) let $(t_k)_k\in \CFE$. We suppose that $t_k=\infty$ for $k>k_0$ and $t_{k_0}=1$. Then :
   % $$ \forall n\geqslant k_0, [t_0,\cdots,t_n]= [t_0,\cdots,t_{k_0}]$$
   % \esp  since by obvious induction on $n$, $[t_0,\cdots,t_n,\infty]=[t_0,\cdots,t_n]$ for all $n\in\N$. So, this sequence $([t_0,\cdots,t_n])_n$ converges toward $\theta=[t_0,\cdots,t_{k_0}]$.      
    Consequence of the definition and the well known fact : 
    $ \forall m\in \N, [a_0,\cdots,a_{s-1},m]=\frac{mp_{s-1}+p_{s-2}}{mq_{s-1}+q_{s-2}}$.    
     \end{mademo}

      $\bullet$ Let $\alpha$ be a rational and $[a_0,a_1,a_2,\cdots,a_r,1]$ its CFE. ( we denote $a_{r+1}=1$)\\
      \esp We have the following induction formula :
       $$ p_{-2}=0 \esp ; \esp  p_{-1}=1 \esp ; \esp  \forall n\in \inter{0}{r+1} \esp , \esp p_n=a_np_{n-1}+p_{n-2}$$
       $$ q_{-2}=1 \esp ; \esp q_{-1}=0 \esp ; \esp  \forall n\in \inter{0}{r+1} \esp , \esp q_n=a_nq_{n-1}+q_{n-2}$$
       
     \esp   We have  $\alpha=\frac{p_r+p_{r-1}}{q_r+q_{r-1}}=\frac{p_{r+1}}{q_{r+1}}$. \\
   %    \esp We denote for all $n\in \inter{0}{r}$, $m_n=q_n+q_{n-1}$, for this will be important for our system of numeration. \\
       
       \esp Let $\alpha'=[a'_0,a'_1,\cdots,a'_{r'},1]$ be an other rational with $r'\geqslant r$. With obvious notations, we see that , for  $n\in\inter{0}{r}$ :
       $$(\forall k\in\inter{0}{n}, a_k\leqslant a'_k) \Rightarrow \left(\forall k\in\inter{0}{n}, p_k\leqslant p'_k \text{ and } q_k\leqslant q'_k \right)$$
        \esp In addition, for $j,n$ integers such that $1\leqslant j\leqslant n\leqslant r$ :
        $$(a_j<a'_j \text{ and } \forall k\in\inter{1}{j-1}, a_k\leqslant a'_k) \Rightarrow q_j<q'_j$$
  %   \esp We remark with no surprise, that $a_0$ has no influence on the denominators $(q_i)_i$.\\
        
\newpage

  $\bullet$ Now, we would like to precise the CFE of reals in $\overset{\longleftrightarrow}{[\theta, \theta']}$ ( denotes the set of reals that are between $\theta$ and $\theta'$, even if $\theta > \theta'$), where $\theta$ and $\theta'$ are two different reals and find the rationals in this interval with the lowest reduced denominator.    \\

   \esp First, we introduce a simple and natural notion :
   
   \begin{madef} [ CFE-depth of a real]. \\
   let $x$ be a real. We name \emph{ CFE-depth} of $x$ the non negative integer, denoted $\mu(x)$ and defined by :   $\mu(x)=+ \infty$ if $x$ is  irrational  and $\mu(x)=s$, if $x=[a_0,a_1,\cdots,a_s,1]$ is the CFE of $x$.\end{madef} 
   
   \esp We  remark that : 
   $$ \mu(x)=0 \Leftrightarrow x\in\Z \esp ; \esp \forall n\in\Z, \mu(x+n)=\mu(x) \esp ; \esp \forall x\not\in \Z, \mu(T(x))=\mu(x)-1$$
   
    \esp We denote $\theta = [t_k]_{k\in\N}$ and       $\theta' = [t'_k]_{k\in\N}$, according to our $\phi_1$-representation. We will abreviate $t$ and $t'$ these CFE-sequences. We denote $r$ the smallest integer $k$ such that $t_k\not = t'_k$. Then we have $r\leqslant \min(\mu(\theta),\mu(\theta'))+2$, when $\theta$ or $\theta'$ is rational ( if they are both irrationals, $r$ is finite ! ). Indeed, the extremal case when $r=\mu(\theta)+2$ for example corresponds to $\theta=[t_0,\cdots,t_{r-2},1]$ and $\theta'=[t_0,\cdots,t_{r-2},1,t'_r,...]$, with $t'_r<\infty$.

%  $\bullet$ We now deal  with the " convergents point of view " :  suppose that $\theta'<\theta$. Then, these two reals have convergent sequences $(p_k/q_k)_k$ and $(p'_k/q'_k)_k$ respectively, such that ( $\epsilon\in \{ -1,1\}$) :
%  $$ \frac{p_0}{q_0} < \frac{p_2}{q_2}< \cdots < \frac{p_{2j}}{q_{2j}}\leqslant \theta \leqslant \frac{p_{2j+\epsilon}}{q_{2j+\epsilon}}<\cdots \frac{p_3}{q_3}<\frac{p_1}{q_1}$$ 
%  \esp We do not consider $\theta$ as a proper convergent of itself, if $\theta$ is rational. \\
 % \esp We have the same thing for $\theta'$, with obvious notations...\\
 % \esp These two convergent-sequences of $\theta$ and $\theta'$ coincides until a certain index $r-1$, so, with $i=\lfloor \frac{r-1}{2}\rfloor$ and $\epsilon =1$ if $r$ is even, and $-1$ else :
 % $$ \frac{p_0}{q_0} < \frac{p_2}{q_2}< \cdots < \frac{p_{2i}}{q_{2i}}\leqslant \theta'< \theta \leqslant \frac{p_{2i+\epsilon}}{p_{2i+\epsilon}}<\cdots \frac{p_3}{q_3}<\frac{p_1}{q_1}$$
  
 % \esp Remark that if $k\leqslant \mu(\theta)+1, \mu(\theta')+1$, then : 
 % $$ \forall j\in \inter{0}{k}, \frac{p_j}{q_j}=\frac{p'_j}{q'_j} \Leftrightarrow  \forall j\in \inter{0}{k},t_j=t'_j$$
 % \esp Indeed :  $(\Leftarrow)$ is obvious by induction. $(\Rightarrow)$ is deduced from the fact that the convergents are reduced, so : forall $j\in\inter{0}{k}$, $p_j=p'_j$ and $q_j=q'_j$, so $t_j=t'_j$ by obvious induction. \\
  %\esp With " proper convergent" of a rational $\theta$, we will mean a convergent of $\theta$, different from $\theta$.\\
   \esp We remark that, all integers in $\overset{\longleftrightarrow}{[\theta, \theta']}$ minimize the denominator of their reduced fraction : it is $1$ !! So, we can suppose that $\lfloor \theta \rfloor = \lfloor \theta' \rfloor$ and even that $\theta, \theta'\in [0,1[$.\\
   \esp The following Lemma proves that, in that case, there is only one rational in $\overset{\longleftrightarrow}{[\theta, \theta']}$, that minimizes the value of its denominator : it is usually named  the " best rational" in $\overset{\longleftrightarrow}{[\theta, \theta']}$ 
   
   \begin{maprop} let $\theta$ and $\theta'$ be two different reals in $[0,1[$ and       $\theta = [t_k]_{k\in\N},\theta' = [t'_k]_{k\in\N}$  their respective CFE. We denote $r$ the lowest integer $k$ such that $t_k\not = t'_k$. \\
   \textbf{(i)} there is a unique rational in  $\overset{\longleftrightarrow}{[\theta, \theta']}$ that minimizes the denominator. We denote it $\gamma$.\\
   - if $r\leqslant \min(\mu(\theta),\mu(\theta'))$, then $\gamma= [t_0,\cdots,t_{r-1},\min(t_r,t'_r),1]$.\\
   - else, $\mu(\theta)<\mu(\theta')$ ( up to swap) and $\gamma=\theta$.\\
  \textbf{(ii)} in both cases, $\mu(\gamma)\leqslant \min(\mu(\theta),\mu(\theta'))$ and $\gamma = [t_0,\cdots,t_{s-1},\min(t_s,t'_s),1]$, where $s=\mu(\gamma)\leqslant r$  and $\forall k\in\inter{0}{s-1},t_k=t'_k$.\\
 \textbf{(iii)} the best rational in  $\overset{\longleftrightarrow}{[\theta, \theta']}$ is the common semi-convergent of $\theta$ and $\theta'$ with the greatest denominator.
   
 %  (ii) convergents point of view :\\
   % let $k\leqslant \mu(\theta)+1, \mu(\theta')+1$ and  $\frac{p_j}{q_j}$ and $\frac{p'_j}{q'_j}$ be the $j^{th}$ convergent of $\theta$ and $\theta'$ ( resp. ). Then : 
  % $$ \forall j\in \inter{0}{k}, \frac{p_j}{q_j}=\frac{p'_j}{q'_j} \Leftrightarrow  \forall j\in \inter{0}{k},t_j=t'_j \Leftrightarrow \forall j\in \inter{0}{k}, \frac{p_j}{q_j},\frac{p'_j}{q'_j} \not\in ]\theta,\theta'[$$
 % - if every proper convergent of $\theta$  is a convergent of $\theta'$ then $\theta$ is the best rational in $\overset{\longleftrightarrow}{[\theta, \theta']}$ ( and vice versa).\\
%- else, let $r$ be the lowest integer $j$ such that $\frac{p_j}{q_j}\not =\frac{p'_j}{q'_j}$. Then, only one of them is in $\overset{\longleftrightarrow}{[\theta, \theta']}$. If we denote $p"_r/q"_r$ those that is not in $\overset{\longleftrightarrow}{[\theta, \theta']}$, then the best rational in $\overset{\longleftrightarrow}{[\theta, \theta']}$ is $\frac{p"_r+p_{r-1}}{q"_r+q_{r-1}}$.   
  
   \end{maprop}    
%\textbf{ Remark} : we can deduce from this Lemma, that the best rational in  $\overset{\longleftrightarrow}{[\theta, \theta']}$ is  $\theta$ if and only if (($r>\mu(\theta)$ and $\mu(\theta')>\mu(\theta)$ ) or $( r=\mu(\theta)$ and $t_r<t'_r)$).

   \begin{mademo}
   (i)  if $r\leqslant \min(\mu(\theta),\mu(\theta'))$. Suppose that $t_r<t'_r$. We have for $(d_k)_{k\in\N}\in \CFE$:
     $$ [d_k]_{k\in\N}\in\overset{\longleftrightarrow}{[\theta, \theta']} \Leftrightarrow \begin{cases}  \forall k<j, d_k=t_k=t'_k  \\ \sigma^{r}(t) \leqslant_A  \sigma^{r}(d)\leqslant_A \sigma^r(t') \esp (*)\end{cases}$$ 
      \esp where $\sigma$ is the usual shift : for any sequence $u$, $\forall k\in\N,\sigma(u)_k=u_{k+1}$.\\
   \esp But, if we want the lowest denominator for the rational $[d_k]_{k\in\N}$, we have to choose the lowest $d_k$ or the $\infty$ value ( if possible), for all $k$. So we have to choose first $d_r=t_r$ and then, the condition (*) becomes : $\sigma^{r+1}(d)\leqslant_A \sigma^{r+1}(t)$. So, we choose $d_{r+1}=1$ and $\forall k>r+1,d_k=\infty$. \\
    - else, one at least of $\mu(\theta)$ and $\mu(\theta')$ is finite and they can not be equal, since $r$ can not be greater than both of them. Suppose $\mu(\theta)< \mu(\theta')$, then we have $\mu(\theta)<r$ and $\forall k\in\inter{0}{\mu(\theta)}, t_k=t'_k$. So, the same arguments as in the previous case prove that $\theta$ is the best rational in $\overset{\longleftrightarrow}{[\theta, \theta']}$.\\
   
   (ii)  it is plain in the first case, since $\mu(\gamma)=r$. If $\mu(\theta)<r$ and $\mu(\theta)<\mu(\theta')$, then $\gamma=\theta$ and $t_s=t'_s$. \\
   
   (iii) is a consequence of (ii), Lemma 1 and the remark following it.
%   (ii)\\
   % The first equivalence is almost plain : $\Leftarrow$ is obvious by induction. $\Rightarrow$ is deduced from the fact that the convergents are reduced, so : forall $j\in\inter{0}{k}$, $p_j=p'_j$ and $q_j=q'_j$, so $t_j=t'_j$ by obvious induction. \\
  % \esp The second equivalence : $(\Rightarrow)$ : if $t_j=t'_j$ for all $j\in\inter{0}{k}$, then :
 %  $$ \forall j\in\inter{0}{k}, \frac{p_j}{q_j}=\frac{p'_j}{q'_j}=[t_0,\cdots,t_j,\overline{\infty}]$$
  % \esp This is a $\phi$-CFE if $t_j>1$ and a $\phi_1$-CFE if $t_j=1$. Let denote $\alpha_j$ the former rational. We have 2 cases : \\
  % --- Case 1 : $r$ is even, then $\alpha_j$ 
%   - if every proper convergent of $\theta$ is a convergent of $\theta'$ then $r>\mu(\theta)$ and $\mu(\theta')>\mu(\theta)$, with notations of (i), so we are in the second case of (i). \\
%   - else we have $t_j=t'_j$ for $j\in\inter{0}{r-1}$ and $t_r\not = t'_r$ ( see the remark above). So $r$ is as in (i) . In addition, $r\leqslant \min(\mu(\theta),\mu(\theta'))$ since at least one proper convergent of $\theta$ is not a convergent of $\theta'$. So, the "best rational" in $\overset{\longleftrightarrow}{[\theta, \theta']}$ is  $\alpha=[t_0,\cdots,t_{r-1},\min(t_r,t'_r),1]$. With usual notations, we have :
 %  $$ \alpha = \frac{p_{r-1}+p"_r}{q_{r-1}+q"_r}$$
  % \esp with : $p"_r=p_r, q"_r=q_r$ if $t_r<t'_r$ or $p"_r=p'_r, q"_r=q'_r$, else. %.................... 
   \end{mademo}
  
  \textbf{Remark : } as a direct consequence of (iii) : $\theta$ is the best rational in $\overset{\longleftrightarrow}{[\theta, \theta']} $ if and only if $\theta$ is a semi-convergent of $\theta'$. \\

  $\bullet$  Let $\alpha$ be a real, $[a_k]_{k\in \N^*}$ its CFE and $r=\mu(\alpha)$, the CFE-depth of $\alpha$. So, we denote $[a_0,a_1,\cdots,a_r,1]$ the CFE of $\alpha$ if $\alpha$ is rational. We also denote $(p_n/q_n)_n$ the usual sequence of convergents of $\alpha$.\\
 \esp We consider the usual notion of best rational approximation of a real $\alpha$ : for $p,q$ two integers, $p/q$ is said a \emph{best rational approximation} of $\alpha$ if and only if :
 $$ \forall q'\in \inter{1}{q-1}\esp , \esp \forall p'\in\Z \esp , \esp  \left|\dfrac{p'}{q'}-\alpha \right| > \left|\dfrac{p}{q}-\alpha \right|$$
 \esp It is well known that best rational approximation of a real are exactly its reduced convergents. \\
 \esp Now, we can consider two sided similar definitions : 
 
 \begin{madef} [ best sided rational approximation] .\\
 for $p,q$ two integers, $p/q$ is said a \emph{best left rational approximation} of $\alpha$ if and only if :
 $$ \forall q'\in \inter{1}{q-1}\esp , \esp \forall p'\in\Z \esp , \esp  \dfrac{p'}{q'}< \dfrac{p}{q} \leqslant \alpha \esp \text{ or } \esp  \dfrac{p'}{q'}>\alpha$$
  \esp $p/q$ is said a \emph{best right rational approximation} of $\alpha$ if and only if :
 $$ \forall q'\in \inter{1}{q-1}\esp , \esp \forall p'\in\Z \esp , \esp  \dfrac{p'}{q'} > \dfrac{p}{q} \geqslant \alpha \esp \text{ or } \esp  \dfrac{p'}{q'}<\alpha$$
 \end{madef}

 \esp Here is a corollary of Proposition 1 :
       
  \begin{moncoro} .\\
  (i) best left rational approximations of $\alpha$ are the semi-convergents of $\alpha$, that are lower than $\alpha$.\\
   (ii) best right rational approximations of $\alpha$ are the semi-convergents of $\alpha$, that are greater than $\alpha$. 
   \end{moncoro}
   \begin{mademo}(i) we remark that $p/q$ is a best left rational approximation of $\alpha$ if and only if $p/q$ is the best rational in $[p/q, \alpha]$ and use the remark below Proposition 1. Same arguments for (ii).   
   \end{mademo}    
       
    \esp If  we denote  $(p_k/q_k)_k$ the reduced convergents of $\alpha$, then : \\
    - its best left rational approximations are  : 
    $$\frac{p_{2i}+mp_{2i+1}}{q_{2i}+mq_{2i+1}} \esp ; \esp i\in \inter{0}{(\mu(\alpha)-1)/2} \esp ; \esp m\in \inter{0}{a_{2i+2}}$$
     - its best right rational approximations are  : 
    $$\frac{p_{2i-1}+mp_{2i}}{q_{2i-1}+mq_{2i}} \esp ; \esp i\in \inter{1}{\mu(\alpha)/2} \esp ; \esp m\in \inter{0}{a_{2i+1}}$$

%  $\bullet$ The sequence $(\{ n\alpha \})_n$ has been widely studied ( see .....) for irrational and for rational $\alpha$. \\
 % \esp The links with continued fraction expansion of $\alpha$ and Ostrowski numeration is well known ( see ....). For convenience, we will adopt a variation of these numeration system : \\

 \newpage

\section{A numeration system}

    \subsection{Ostrowski's numeration}
 \label{sec:ost}    
    
 \esp We will only deal here with the case $\alpha$ irrational, even if the rational case is interesting ( see next section). We denote $\Omega_{\alpha}$ the set of sequences of integers defined as follows ( we denote $[a_k]_{k\in\N}$ the continued fraction expansion of $\alpha$) :
   $$ \Omega_{\alpha}=\{ (d_n)_{n\in\N^*}, d_1\in \inter{0}{a_1-1}, \forall k\in\N^*\backslash \{ 1 \}, d_k\in \inter{0}{a_k} \text{ and } ( d_k=a_k\Rightarrow d_{k-1}=0 ) \}$$
   
  \esp What we call " markovian condition"  is the last implication :    $d_k=a_k\Rightarrow d_{k-1}=0$.\\
  \esp From this set of infinite sequences, we extract two subsets, that will be our numeration sets for reals and integers respectively : 
  $ O_{\alpha}$ is the set of sequences $d$ of $\Omega_{\alpha}$ such that $d$ does not " end with" $(\max,0)^{\infty}$, an infinite sequence $a_k0a_{k+2}0\cdots$. So to say, there is an infinite number of even and an infinite number of odd values of $k$ such that $d_k<a_k$. Now, $ O_{(\alpha)}$ is the set of sequences $d$ of $\Omega_{\alpha}$ ( or $O_{\alpha}$) that ends with an infinite sequence of $0$ : so to say $d_k=0$ for any sufficiently large $k$.\\
  
  \esp We define then two maps :
  $$ f_{\alpha} : \begin{cases} O_{(\alpha)}\to \N \\ 
                                d \to \sum\limits_{k=1}^{\infty}d_kq_{k-1}\end{cases}  \esp ; \esp  
  g_{\alpha} : \begin{cases} O_{\alpha}\to [-\alpha,1-\alpha[ \\ 
                                d \to \sum\limits_{k=1}^{\infty}d_k(\alpha q_{k-1}-p_{k-1}) \end{cases} $$

    \esp It is well known that $f_{\alpha}$ and $g_{\alpha}$ are well defined and are bijective. Moreover :
    $$ \forall d\in O_{(\alpha)}, \esp \{ f_{\alpha}(d) \alpha \}=\{ g_{\alpha}(d) \} $$
    \esp But, we will emphasize an other aspect : the maps above are increasing for the usual order on $\N$ and $\R$ respectively and following orders on $O_{(\alpha)}$ and $O_{\alpha}$. \\

      - the \textbf{reversed lexicographic order ( RLO) } on  $O_{(\alpha)}$ : 
     $$ d \underset{R}{\leqslant} d' \Leftrightarrow d=d' \text{ or } \exists j\in \N^*, \begin{cases} d_j<d'_j \\
     \forall i>j, d_i=d'_i \end{cases} $$
      - the \textbf{alternate lexicographic order ( ALO)} on $ O_{\alpha}$ : 
     $$ d \underset{A}{\leqslant}d' \Leftrightarrow d=d' \text{ or } \exists j\in \N^*, \begin{cases} (-1)^{j-1}d_j<(-1^{j-1}d'_j \\
     \forall i\in \inter{1}{j-1}, d_i=d'_i \end{cases} $$
        
   \esp These are total orders on these sets respectively. \\

   \esp Our aim is to find a variant of Ostrowski numeration that has same properties, but that code reals of $[0,1[$ instead of $[-\alpha, 1-\alpha[$ and also all integers, positive as well as negative ones.\\ \esp We will see that it suffices to change the markovian condition : instead of  $d_k=a_k\Rightarrow d_{k-1}=0$, we take $d_k=0 \Rightarrow ( d_{k-1}=a_{k-1}$ or $d_i=0$ for all $i\geqslant k$).
\newpage

      \subsection{$\alpha$-numeration for a rational $\alpha$}
      
\label{sec:rational-numeration}   
       
   \esp Why do we consider this case $\alpha$ rational ? Indeed,  the set $\{ \{ n\alpha \}, n\in\N \}$ is finite and trivial. It  can not define a base of numeration for $[0,1[$. But the order properties of the sequence $(\{ n\alpha \})_{n\in\N}$ are not obvious and our Ostrowski-like numeration will help. \\
 %  \esp Moreover, this will be very useful in two other articles about some numerical semigroups ( see ....).\\  

   $\bullet$ Let $\alpha$ be a rational in $[0,1[$ and $\alpha=[0,a_1,\cdots,a_r,1]$ its CFE. We will denote $(p_k/q_k)_{0\leqslant k \leqslant r+1 }$ its convergents, so that $\alpha=\frac{p_{r+1}}{q_{r+1}}$.
   
   \begin{madef} [ $\alpha$-admissible sequences].\\
    a sequence $d$ in $\N^r$ is said \emph{$\alpha$-admissible} if and only if :
   $$\forall j\in\inter{1}{r},\begin{cases}  d_j\in\inter{0}{a_j}\\  d_j=0 \Rightarrow ( \forall i\geqslant j, d_i=0 )\text{ or } d_{j-1}=a_{j-1} \end{cases}$$
   \end{madef}
   
    \esp We will denote $E_{\alpha}$ the set of $\alpha$-admissible sequences. \\
   
   \textbf{Remark : } for $j=1$, the second condition reduces to $d_1=0\Rightarrow \forall i\geqslant 1,d_i=0$. So to say, $d=(0,\cdots,0)$ is the only element of $E_{\alpha}$ , whose first coordinate is $0$.

     \begin{monlem}.\\
     (i) 
     $$ \forall d\in E_{\alpha}, \forall k\in\inter{1}{r}, \esp \sum_{i=1}^kd_iq_{i-1}< q_k+q_{k-1}$$
     
     (ii) let $d,d'\in E_{\alpha}$ and $n\in\inter{1}{r}$ such that $d'_n>0$. 
     $$ \forall k\in \inter{1}{n}, \esp
     \sum_{i=1}^k (d_i-d'_i)q_{i-1} < q_k $$
     \end{monlem}
     \begin{mademo}
     (i) by plain induction on $k$. \\
     (ii) by induction ( on 2 ranks) on $k$ : \\
     - it is true for $k=0$ ( obvious) and for $k=1$ : indeed $d'_1>0$ ( else, we would have $d'=0$ and $d'_n=0$) : then $(d_1-d'_1)q_0\leqslant (a_1-1)q_0 = q_1-1$. \\
     - we suppose that it is true for the ranks $k-2$ and $k-1$, where $k$ is an integer in $\inter{2}{n}$. Then, we have two cases for the rank $k$ : \\
     $\blacktriangleright$ Case 1 : $d_k-d'_k\leqslant a_k-1$, then, with the induction hypothesis on rank $k-1$ :
     $$   \sum_{i=1}^k (d_i-d'_i)q_{i-1}< q_{k-1}+(a_k-1)q_{k-1}=a_kq_{k-1} = q_k-q_{k-2} < q_k $$
     \esp the last inequality is true, for $k\geqslant 2$. \\
      $\blacktriangleright$ Case 2 : $d_k=a_k, d'_k=0$, then $d'_{k-1}=a_{k-1}$ ( else, we would have $d'_j=0$ for all $j\geqslant k$, but $d'_n\not = 0$) and $d_{k-1}-d'_{k-1}\leqslant 0$. So with the induction hypothesis on rank $k-2$ :
     $$  \sum_{i=1}^k (d_i-d'_i)q_{i-1}< q_{k-2}+a_kq_{k-1}=q_k$$
     \end{mademo}
 
    \esp We consider  the \textbf{reversed lexicographic order ( RLO) } denoted $ \underset{R}{\leqslant}$ on $\N^r$ : 
     $$ d \underset{R}{\leqslant} d' \Leftrightarrow d=d' \text{ or } \exists j\in \inter{1}{r}, \begin{cases} d_j<d'_j \\
     \forall i\in \inter{j+1}{r}, d_i=d'_i \end{cases} $$
    \esp  It is a  total order on $E_{\alpha}$.
       
     \begin{monlem}
       the map $\Psi_{\alpha}$ below is increasing from $(E_{\alpha} ,\leqslant_R)$ to $(\inter{0}{q_{r+1}-1},\leqslant)$.
     $$ \Psi_{\alpha} : \begin{cases} E_{\alpha} \to \inter{0}{q_{r+1}-1} \\ d\to \sum\limits_{j=1}^rd_jq_{j-1} \end{cases}$$
      \end{monlem}
      
      \begin{mademo}
      First,  for all $d\in E_{\alpha}, \Psi_{\alpha}(d) \in \inter{0}{q_{r+1}-1}$, with Lemma 2 (i). \\
   \esp  Now, let prove that $\Psi_{\alpha}$ is increasing. Let $d,d'\in E_{\alpha}$, such that $d<_R d'$. We have $j\in \inter{1}{r}$, such that :
   $$d_j<d'_j \esp \text{ and } \esp \forall i\in\inter{j+1}{r}, d_i=d'_i$$
   \esp So : 
  $$ \Psi_{\alpha}(d')-\Psi_{\alpha}(d)= \sum_{i=1}^{j-1}(d'_i-d_i)q_{i-1}+ (d'_j-d_j)q_{j-1} $$
  \esp We just have to prove that : $\sum\limits_{i=1}^{j-1}(d_i-d'_i)q_{i-1}<q_{j-1}$, since $d'_j-d_j\geqslant 1$. This is shown by  Lemma 2 (ii), for $d'_j>0$.
      \end{mademo}

      \esp Now, we prove that $\Psi_{\alpha}$ is surjective :  the following algorithm explains the inverse function of $\Psi_{\alpha}$. We will denote $m_k=q_k+q_{k-1}$ for any $k\in\inter{0}{r}$. So $m_r=q_{r+1}$.

  \begin{monalgo}  let $n\in \inter{0}{m_r-1}$. \\
  With the following algorithm, we have  $d\in E_{\alpha}$ and $\Psi_{\alpha}(d)=n$.
   \end{monalgo}
    \begin{algorithm}
%\caption{Voici un algorithme.}
    \KwIn{$n$}
    \KwOut{$(d_i)_{i\in \inter{1}{r}}$}
\BlankLine

\For{$k \gets r$ \KwTo $1$ \KwStep $-1$}{
 $d_k \gets \max\left(0,\left\lfloor \frac{n-q_{k-2}}{q_{k-1}}\right\rfloor \right)$ \; 
 $ n \gets n-d_kq_{k-1}$}
%\Return{$x$}
\end{algorithm}
 
   \begin{mademo}
   \esp We begin with a remark : if $n<m_s$ for an integer $s\in\inter{1}{r}$, then : $d_k=0$ for $k\in\inter{s+1}{r}$. Indeed, we will have $n<m_k$ for all $k\in \inter{s}{r}$ so $n-q_{k-2}<q_{k-1}$ for all $k\in \inter{s+1}{r}$. \\
   \esp Let us prove the result  by induction on $s$, where $s$ in an integer such that $n\in\inter{0}{m_s-1}$ : \\
     - for $s=1$, $m_1=a_1+1$. Let $n\in\inter{0}{a_1}$. Then $d_1=n$ and $d=(d_1)\in E_{\alpha}, \Psi_{\alpha}(d)=d_1=n$.\\
     - we suppose that the algorithm is available for all $n\in \inter{0}{m_{s-1}-1}$, with $s\geqslant 2$.\\
     \esp  Let $n\in \inter{m_{s-1}}{m_s-1}$.  Then $q_{s-1}\leqslant n-q_{s-2}<(a_s+1)q_{s-1}$, so $ d_s\in\inter{1}{a_s}$. We denote $n_1=n-d_sq_{s-1}$, the value of $n$ after the loop for $k=s$. We have : $ q_{s-2} \leqslant n_1< q_{s-1}+q_{s-2}=m_{s-1}$. By induction hypothesis, $d'=(d_1,\cdots,d_{s-1})\in E_{\alpha}$ and $n_1=\Psi_{\alpha}(d')=\sum\limits_{i=1}^{s-1}d_iq_{i-1}$. But, $n=n_1+d_sq_{s-1}$ and so $\Psi_{\alpha}(d)=n$, because we have $d\in E_{\alpha}$ : indeed, we have 2 subcases : \\
    $\blacktriangleright$ Case 1 : if $n_1\geqslant m_{s-2}$, then $d_{s-1}>0$ and, since $(d_1,\cdots,d_{s-1})\in E_{\alpha}$, then $d\in E_{\alpha}$. \\
     $\blacktriangleright$ Case 2 : if $ n_1< m_{s-2}$ ( which leads to $s\geqslant 3$, for $m_0=q_0$), then $d_{s-1}=0$ and $n_2=n_1$ ( $n_2$ : the value of $n$ after the loop for $k=s-1$). But, since $n_1\geqslant q_{s-2}$, then $n_2-q_{s-4}\geqslant a_{s-2}q_{s-3}$ and finally $d_{s-2}=a_{s-2}$. By induction hypothesis, $d'=(d_1,\cdots,d_{s-3},a_{s-2},0)$ is $\alpha$-admissible, so $d\in E_{\alpha}$. 
   \end{mademo}

  \begin{maprop}
  $\Psi_{\alpha}$ is an order isomorphism between  $(E_{\alpha} ,\leqslant_R)$ and $(\inter{0}{q_{r+1}-1},\leqslant)$.
  \end{maprop}
  
  \textbf{Remark : } as a direct consequence : $E_{\alpha}$ has $q_{r+1}$ elements.
  
  \begin{mademo} a direct consequence of Lemma 3 and Algorithm 1 
  \end{mademo}

     $\bullet$ Now, we will deal with $\alpha$-numeration for elements of  $U_{\alpha}=  \{ \{ k\alpha\}, k\in\N \} $. Since, $\alpha=\frac{p_{r+1}}{q_{r+1}}$ and this fraction is reduced, we have $U_{\alpha}= \{ \frac{n}{q_{r+1}}, n\in \inter{0}{q_{r+1}-1}\}$. So, this set is very simple, but we will focus on the map $k\to \{ k\alpha \}$, with the order point of view :\\
     \esp We consider the \textbf{alternate lexicographic order ( ALO)} denoted $ \underset{A}{\leqslant}$ on $\R^r$ : 
     $$ d \underset{A}{\leqslant}d' \Leftrightarrow d=d' \text{ or } \exists j\in \inter{1}{r}, \begin{cases} (-1)^{j-1}d_j<(-1)^{j-1}d'_j \\
     \forall i\in \inter{1}{j-1}, d_i=d'_i \end{cases} $$
     \esp It is another total order on $E_{\alpha}$. We define also  : 
      $$\forall i\in\inter{-2}{r}, \esp \delta_i=(-1)^i( q_{i}\alpha - p_i)$$
      \esp We have, with $a_0=0$ here :
      $$ \delta_{-2}=\alpha \esp ; \esp \delta_{-1}=1 \esp ; \esp \delta_0=\{ \alpha \} =\alpha \esp ; \esp \forall i\in\inter{0}{r}, \esp \delta_i=-a_i\delta_{i-1}+\delta_{i-2}$$
      \esp Let $T$ be the Gauss map : $]0,1[\to [0,1[, x\to \{ 1/x \}$.\\
      \esp  By induction on $i$, with the fact that : $a_i=\left\lfloor \frac{1}{T^{i-1}(\alpha)} \right\rfloor$ if $i\leqslant r-1$, we obtain :
      $$ \forall i\in\inter{0}{r-1}, \frac{\delta_i}{\delta_{i-1}}=T^i(\alpha) $$
      
      \esp  Beware : for $i=r$, $T^{r-1}(\alpha)=[0,a_r,1]=\frac{1}{a_r+1}$, so :
      $$ \frac{\delta_r}{\delta_{r-1}}= \frac{\delta_{r-2}-a_r\delta_{r-1}}{\delta_{r-1}}= \frac{1}{T^{r-1}(\alpha)}-a_r=1 $$
      \esp So : $\delta_r=\delta_{r-1}$. We will prove ( proof of Algorithm 2) that 
      $ \delta_r=\delta_{r-1}=\frac{1}{q_{r+1}}$.\\
      \esp To summarize this :
      $$ \forall i\in \inter{0}{r-1}, 0<\delta_i<\delta_{i-1} \esp ; \esp \delta_r=\delta_{r-1}=\frac{1}{q_{r+1}}$$
  
     \begin{monlem} let $d,d'\in E_{\alpha}$ and  $j\in\inter{1}{r}$, , then : 
     $$(-1)^{j-1}(d'_j-d_j)>0 \Rightarrow \sum\limits_{i=j+1}^r (-1)^{i}(d'_i-d_i)\delta_{i-1}<\delta_{j-1}$$ 
     \end{monlem}
     
     \begin{mademo}
     First, we remark that, for all $i$, we have $ (-1)^{i}(d'_i-d_i)\leqslant a_i$, so :
     $$  \sum\limits_{i=j+1}^r (-1)^{i}(d'_i-d_i)\delta_{i-1} \leqslant  \sum\limits_{i=j+1}^r a_i\delta_{i-1} = \sum\limits_{i=j+1}^r (\delta_{i-2}-\delta_i)= \delta_{j-1}+\delta_j-\delta_{r-1}-\delta_r$$
 
  $\blacktriangleright$  Case 1 : if $(-1)^{j+1}(d'_{j+1}-d_{j+1})\leqslant a_{j+1}-1$, then :
   $$ \sum\limits_{i=j+1}^r (-1)^{i}(d'_i-d_i)\delta_{i-1} \leqslant  \delta_{j-1}-\delta_{r-1}-\delta_r< \delta_{j-1}$$
    $\blacktriangleright$ Case 2 : if $(-1)^{j+1}(d'_{j+1}-d_{j+1})=a_{j+1}$.\\ 
     $\blacktriangleright \blacktriangleright$ Subcase 1 : if $j$ is even, $d'_{j+1}=0$ and $d_{j+1}=a_{j+1}$. We can not have $d'_j=a_j$, for, with our hypothesis, $d_j>d'_j$. So $d'_i=0$ for all $i>j$ and :
   $$ \sum\limits_{i=j+1}^r (-1)^{i}(d'_i-d_i)\delta_{i-1}=a_{j+1}\delta_j - \sum_{i=j+2}^{r}  (-1)^{i}d_i\delta_{i-1}\leqslant  \sum_{p=0}^{(r-j-1)/2}  a_{j+2p+1}\delta_{j+2p}= $$
$$  = \sum_{p=0}^{(r-j-1)/2} (\delta_{j+2p-1}-\delta_{j+2p+1})= \delta_{j-1}-\delta_{r'}<\delta_{j-1}$$
\esp \esp with $r'=r$ or $r-1$. \\
 $\blacktriangleright \blacktriangleright$ Subcase 2 : if $j$ is odd, similar arguments lead to the same conclusion ( we swap $d$ and $d'$).    
     \end{mademo}
    
     \begin{maprop} .\\
      (i) the map $\Lambda_{\alpha} $ ( defined below) is an order isomorphism, with ALO on $E_{\alpha}$ : 
       $$ \Lambda_{\alpha} : \begin{cases} E_{\alpha} \to \left\{ \frac{n}{q_{r+1}}, n\in \inter{0}{q_{r+1}-1}\right\} \\ d\to \sum\limits_{j=1}^rd_j(-1)^{j-1}\delta_{j-1} \end{cases}$$
       
     (ii) we have : 
     $$ \forall n\in \inter{0}{q_{r+1}-1}, \{ n\alpha \}= \Lambda_{\alpha}(\Psi_{\alpha}^{-1}(n))$$

 %    (iv)  if $\alpha\in [0,1[$, then for $n\in \inter{0}{m_r-1}$, $  \lfloor (n+1)\alpha\rfloor - \lfloor n\alpha\rfloor =1 \Leftrightarrow n=\Psi_{\alpha}(d)$ with $d_k=a_k$ for all $k<j$ and $d_j<a_j$, with $j$ odd.  
     \end{maprop}

  \begin{mademo} 
  (i) First, we will show that $\Lambda_{\alpha}$ is increasing : let $d,d'\in E_{\alpha}$ with $d<_A d'$. Then, we have $j\in\inter{1}{r}$ such that :
  $$(-1)^{j-1}d_j<(-1)^{j-1}d'_j \esp \text{ and } \esp \forall i<j, d_i=d'_i$$
  \esp  So :
  $$ \Lambda_{\alpha}(d')-\Lambda_{\alpha}(d)= (-1)^{j-1}(d'_j-d_j)\delta_{j-1} +\sum_{i=j+1}^r (-1)^{i-1}(d'_i-d_i)\delta_{i-1} $$
  \esp Now :
  $$(-1)^{j-1}(d'_j-d_j)\delta_{j-1}\geqslant \delta_{j-1}$$
  \esp so with Lemma 4, we obtain :
  $$ \Lambda_{\alpha}(d')-\Lambda_{\alpha}(d)>0$$
\esp Now that we have proved that $\Lambda_{\alpha}$ is increasing, we can easily deduce that $\Lambda_{\alpha}(E_{\alpha})\subset [0,1[$ : first, remark that $(0,\cdots,0)$ is the lowest element of $E_{\alpha}$ ( with ALO), so $\Lambda_{\alpha}(d)\geqslant 0$ for all $d\in E_{\alpha}$. Now, 
$(a_1,0,a_3,0,\cdots)$ is the greatest element of $E_{\alpha}$ for ALO, so : 
$$ \forall d\in E_{\alpha}, \Lambda_{\alpha}(d)\leqslant \sum_{p=0}^{(r-1)/2} a_{2p+1}\delta_{2p}=\sum_{p=0}^{(r-1)/2}(\delta_{2p-1}-\delta_{2p+1})=\delta_{-1}-\delta_{r'}<1$$
\esp  with $r'=r$ or $r-1$. \\

(ii) we just have to show this equality to complete the proof : let $d\in E_{\alpha}$. It is sufficient to prove that $\Lambda_{\alpha}(d)=\{\Psi_{\alpha}(d)\alpha \}$. Now :
$$ \Lambda_{\alpha}(d)= \sum\limits_{j=1}^rd_j (q_{j-1}\alpha - p_{j-1})= \alpha \Psi_{\alpha}(d)- k$$
\esp where $k= \sum\limits_{j=1}^rd_j p_{j-1}$ is an integer. So $\Lambda_{\alpha}(d)=\{\Psi_{\alpha}(d)\alpha \}$ modulo 1. But, we have seen that both terms are in $[0,1[$, q.e.d.
%(iv) let $n\in \inter{0}{m_r-1}$. We denote $n=\Psi_{\alpha}(d)$ and $n+1=\Psi_{\alpha}(d')$. \\
%\esp If  $d_1<a_1$, then  $d'_k=d_k$ for all $k>1$, and $d'_1=d_1+1$. So 
%$ \lfloor (n+1)\alpha\rfloor - \lfloor n\alpha\rfloor =0$, for $p_0=a_0=0$.\\
%\esp  If $d_1=a_1$, then we denote $j$ the smallest integer $k$ such that $d_k<a_k$. We have then : $d'_j=d_j+1$ and $d_{j-1}=0, d_{j-2}=a_{j-2}, d_{j-3}=0$ and so on...except for $d_1=1$ if $j$ is even. But, the value of $d_1$ and $d'_1$ has no importance for the result. So that : 
%$$ \lfloor (n+1)\alpha\rfloor - \lfloor n\alpha\rfloor = p_{j-1}- \sum_{i=1}^{\lfloor (j-1)/2\rfloor}a_{j+1-2i}p_{j-2i}=  p_{j-1}- \sum_{i=1}^{\lfloor (j-1)/2\rfloor}(p_{j+1-2i}-p_{j-2i-1}) $$ 
%\esp the previous expression is $p_1=1$ if $j$ is odd or is $p_0$ if $j$ is even.
 \end{mademo}

    \textbf{Remarks : } result (ii) means that the map $n\to \{n\alpha \}$ ( with $0\leqslant n < q_{r+1}$), is, from the order point of view, the " same thing"  as the identity $(E_{\alpha},RLO)\to (E_{\alpha},ALO)$. \\
     
     \esp We can sum up these formulae  : $ \forall n\in \inter{0}{q_{r+1}-1}$, with $d=\Psi_{\alpha}^{-1}(n)$ :
     $$n=\sum_{j=1}^r d_jq_{j-1} \esp ;\esp   \lfloor n\alpha\rfloor =\sum_{j=1}^r d_jp_{j-1}  \esp ; \esp \{ n\alpha \}=\sum_{j=1}^r (-1)^{j-1}d_j\delta_{j-1}$$

  \esp   The following algorithm expresses the inverse function of $\Lambda_{\alpha}$.

  \begin{monalgo} let $\beta\in \{ \frac{n}{q_{r+1}}, n\in\N \}$. Applying the algorithm below, we have :\\
  (i) $b\in E_{\alpha}$.\\
  (ii)  $\beta = \Lambda_{\alpha}(b)$.
   \end{monalgo}
   
      \begin{algorithm}
%\caption{Voici un algorithme.}
    \KwIn{$\beta$}
    \KwOut{$(b_i)_{i\in \inter{1}{r}}$}
\BlankLine
 
\For{$k \gets 1$ \KwTo $r$}{
 $b_k \gets \min\left(a_k,\left\lceil \frac{\beta}{\delta_{k-1}}\right\rceil \right)$ \; 
 $ \beta \gets b_k\delta_{k-1} -\beta$}
%\Return{$x$}
\end{algorithm}

\begin{mademo} First, we denote $(\beta_k)_{k\in\inter{0}{r}}$ the finite sequence defined by :
$$ \beta_0=\beta \esp ; \esp \forall k\in\inter{1}{r} \esp , \esp \beta_k=b_k\delta_{k-1}-\beta_{k-1}$$
\esp Thus, $\beta_k$ is the value of $\beta$ after $k$ loops in Algorithm 2. So, we have : 
$$b_k=\min(a_k,\lceil \beta_{k-1}/\delta_{k-1}\rceil )$$

 (i) let us verify that $b\in E_{\alpha}$ : by induction on $k$, we will prove that "$(b_1,\cdots,b_k)$ is $\alpha$-admissible  and that $-\delta_k<\beta_k < \delta_{k-1}$ for all $k\in\inter{0}{r}$". \\
 - it is true for $k=0$, since $\delta_0=\alpha>0$ and $\delta_{-1}=1$. \\
 - we suppose that it is true for $k-1$ with $k\in \inter{1}{r}$. 
 Then, $ \frac{\beta_{k-1}}{\delta_{k-1}}>-1$, so $\left\lceil \frac{\beta_{k-1}}{\delta_{k-1}}\right\rceil \geqslant 0$ and $0\leqslant b_k\leqslant a_k$. If $b_{k-1}>0$, then $(b_1,\cdots,b_k)$ is $\alpha$-admissible, for $(b_1,\cdots,b_{k-1})\in E_{\alpha}$. If $b_{k-1}=0$, then $\beta_{k-2}\leqslant 0$ and we have 2 cases : \\
 $\blacktriangleright$ Case 1 : if $\beta_{k-2}=0$, then by obvious induction, $\beta_i=0$ and $b_i=0$ for all $i\geqslant k-1$.\\
 $\blacktriangleright$ Case 2 : else, we have $ \beta_{k-2}< 0$ , so $k\geqslant 3$ and :
 $$\beta_{k-3}=b_{k-2}\delta_{k-3}-\beta_{k-2}> b_{k-2} \delta_{k-3}$$
 \esp So $\frac{\beta_{k-3}}{\delta_{k-3} }>b_{k-2}$, which leads to $b_{k-2}=a_{k-2}$. \\
\esp In these both cases, we have :  $(b_1,\cdots,b_k)$ satisfies the conditions of $E_{\alpha}$. \\
\esp Now : by induction hypothesis, we have : 
$$- \delta_{k-1}< \beta_{k-1}< \delta_{k-2} $$
$\blacktriangleright$ Case 1 : if $ \frac{\beta_{k-1}}{\delta_{k-1}}\leqslant a_k$, then : $b_k=\left\lceil \frac{\beta_{k-1}}{\delta_{k-1}}\right\rceil$ so $\beta_{k-1}\leqslant b_k\delta_{k-1} $, so $\beta_k\geqslant 0$ and $b_k< \frac{\beta_{k-1}}{\delta_{k-1}}+1$, so $\beta_k<\delta_{k-1}$. \\
$\blacktriangleright$ Case 2 : if $\frac{\beta_{k-1}}{\delta_{k-1}}>a_k$, then $b_k=a_k$ and $\beta_k<0$. Moreover :
$$\beta_k=a_k\delta_{k-1}-\beta_{k-1}>a_k\delta_{k-1}-\delta_{k-2}=-\delta_{k}$$

(ii) 
$$ \Lambda_{\alpha}(b)=\sum_{k=1}^r(-1)^kb_k\delta_{k-1}= \sum_{k=1}^r((-1)^k\beta_k-(-1)^{k-1}\beta_{k-1})=(-1)^r\beta_r-\beta  $$
\esp we also have $\beta = \Lambda_{\alpha}(b)+(-1)^r\beta_r$. Now, $-\delta_r<\beta_r < \delta_{r-1}$.\\
\esp Claim : " for every $k\in \inter{-1}{r-1}, q_{r+1}\delta_k$ is the $k^{th}$ remainder, denoted $\rho_k$ in the euclidean algorithm between $p_{r+1}$ and $q_{r+1}$ and we have $\rho_{r-1}=1$." \\
\esp Indeed, by double induction on $k$ : \\
- it is true for $k=-1$ and $k=0$, since $ q_{r+1}\delta_{-1}=q_{r+1}=\rho_{-1}$ and $  q_{r+1}\delta_0=p_{r+1}=\rho_0$. \\
- then, both sequences satisfy the same double induction formula : 
$$ \forall k\in\inter{1}{r-1}, q_{r+1}\delta_k= q_{r+1}\delta_{k-2}-a_kq_{r+1}\delta_{k-1} \esp ; \esp \rho_k=\rho_{k-2}-a_k\rho_{k-1} $$
\esp Now, euclidean algorithm stops when we obtain a rest equal to 0, and the former rest is the greatest common divisor of $\rho_{-1}$ and $\rho_0$, namely $1$ here, since the convergent fractions are reduced. So, $\rho_{r-1}=1$. But, we have chosen the continued fraction expansion of $\alpha$, that ends with $1$, so $\rho_{r-2}=a_r+1$, and $q_{r+1}\delta_r=\rho_{r-2}-a_r\rho_{r-1}=1$. \\
\esp We conclude : $\delta_r=\delta_{r-1}=\frac{1}{q_{r+1}}$ and $|\beta_r|\in [0,1/q_{r+1}[$. \\
\esp Now, the former facts show that $\beta_kq_{r+1}\in\Z$ for all $k$, so : $\beta_r=0$. 
\end{mademo}

$\bullet$ We can easily extend this numeration to $[0,1[$, by adding a last " digit" that can range in $[0,1[$. First, we extend the ALO to $E_{\alpha}\times [0,1[$ : $(d,\epsilon)\leqslant_A (d',\epsilon')$ if and only if ($d= d'$ and $\epsilon\leqslant \epsilon'$) or $d<_A d'$.
  
\begin{moncoro}
      the map $\tilde{\Lambda}_{\alpha} $ is an order isomorphism, with ALO on $E_{\alpha}\times [0,1[$ : 
       $$ \tilde{\Lambda}_{\alpha} : \begin{cases} E_{\alpha}\times [0,1[ \to [0,1[ \\ (d,\epsilon)\to \sum\limits_{j=1}^rd_j(-1)^{j-1}\delta_{j-1}+\epsilon \delta_r \end{cases}$$
       \end{moncoro}
    \begin{mademo}
    a direct consequence of Proposition 3.
    \end{mademo}

\textbf{Remark : } if $ \tilde{\Lambda}_{\alpha}(d,\epsilon)=\beta$ then $\epsilon= \{ q_{r+1}\beta\}$, with usual notations. \\

\newpage

 \subsection{$\alpha$-numeration for an irrational $\alpha$}

\label{sec:irrational-numeration}
  
    $\bullet$ Let $\alpha$ be an irrational and $[a_k]_{k\in\N}$ its CFE. We extend our notion of $\alpha$-admissible sequence :

   \begin{madef}[$\alpha$-admissible sequences].\\
    a sequence $d$ in $\N^{\N^*}$ is said \emph{$\alpha$-admissible} if and only if   $ d$  does not end with $(\max,0)^{\infty}$, an infinite sequence of  $ a_k,0,a_{k+2},0,\cdots $ ( so to say there are  an infinite number of even and odd indices $k$ such that $d_k>0$ or $d_{k+1}<a_k$) and :
   $$\forall j\in\N^*,\begin{cases}  d_j\in\inter{0}{a_j}\\  d_j=0 \Rightarrow ( \forall i\geqslant j, d_i=0 )\text{ or } d_{j-1}=a_{j-1} \end{cases}$$
  
   \end{madef}
   \esp Thus, the null-sequence is the only $\alpha$-admissible sequence that begins with $0$.   We denote $ E_{\alpha}$ the set of $\alpha$-admissible sequences and $ E_{(\alpha)}$ the subset of $ E_{\alpha}$ of sequences,  that ends with $0^{\infty}$, an infinite sequence of $0$. \\

    $\bullet$ We consider two lexicographic total order, respectively on $E_{\alpha}$ and $ E_{(\alpha)}$ : \\
     - the \textbf{reversed lexicographic order ( RLO) } on  $E_{(\alpha)}$ : 
     $$ d \underset{R}{\leqslant} d' \Leftrightarrow d=d' \text{ or } \exists j\in \N^*, \begin{cases} d_j<d'_j \\
     \forall i>j, d_i=d'_i \end{cases} $$
      - the \textbf{alternate lexicographic order ( ALO)} on $ E_{\alpha}$ : 
     $$ d \underset{A}{\leqslant}d' \Leftrightarrow d=d' \text{ or } \exists j\in \N^*, \begin{cases} (-1)^{j-1}d_j<(-1^{j-1}d'_j \\
     \forall i\in \inter{1}{j-1}, d_i=d'_i \end{cases} $$
     
     $\bullet$ We define : 
      $$\forall i\in\N\cup \{ -1\} \esp , \esp \delta_i=(-1)^i( q_{i}\alpha - p_i)$$
      with, as usual $p_i/q_i$ being the reduced fraction of the convergent $[a_0,\cdots,a_i]$. We have then :
      $$ \delta_{-1}=1 \esp ; \esp \delta_0=\alpha \esp ; \esp \forall i\in\N^* \esp , \esp  \delta_i=-a_i\delta_{i-1}+\delta_{i-2}$$
      \esp Let $T$ be the Gauss map : $]0,1[\backslash \Q\to ]0,1[\backslash \Q, x\to \{ 1/x \}$.\\
      \esp  By induction on $i$, with the fact that : $a_i=\left\lfloor \frac{1}{T^{i-1}(\alpha)} \right\rfloor$ if $i\in \N^*$, we obtain :
      $$ \forall i\in\N \esp, \esp \frac{\delta_i}{\delta_{i-1}}=T^i(\alpha) $$
     \esp $(\delta_i)_{i\in\N}$ is a decreasing and positive sequence, that converges towards $0$.\\ 
     
  \begin{monlem}let $d,d'\in E_{\alpha}$ and  $j\in\inter{1}{r}$, , then : 
     $$(-1)^{j-1}(d'_j-d_j)>0 \Rightarrow \sum\limits_{i=j+1}^{\infty} (-1)^{i}(d'_i-d_i)\delta_{i-1}<\delta_{j-1}$$
    \end{monlem}
    
  \begin{mademo} \esp We have 2 cases : \\
   $\blacktriangleright$ Case 1 : if $(-1)^{j+1}(d'_{j+1}-d_{j+1})\leqslant a_{j+1}-1$, then :
   $$ \sum\limits_{i=j+1}^{\infty} (-1)^{i}(d'_i-d_i)\delta_{i-1} \leqslant (a_{j+1}-1)\delta_j + \sum\limits_{i=j+2}^{\infty} (-1)^{i}(d'_i-d_i)\delta_{i-1} $$
   \esp  But, nor $d$ nor $d'$ ends with $(\max,0)^{\infty}$, an infinite sequence of " $(a_k,0)$", so :
  $$ \exists k>j+1, (-1)^k(d'_k-d_k)<a_k \esp ; \esp \forall i>j,(-1)^i(d'_i-d_i)\leqslant a_i$$
 
  \esp We deduce :
  $$   \sum\limits_{i=j+2}^{\infty} (-1)^{i}(d'_i-d_i)\delta_{i-1}<  \sum\limits_{i=j+2}^{\infty} a_i\delta_{i-1}=\sum\limits_{i=j+2}^{\infty} (\delta_{i-2}-\delta_i)=\delta_{j+1}+\delta_j$$
  \esp We conclude : 
  $$\sum\limits_{i=j+1}^{\infty} (-1)^{i}(d'_i-d_i)\delta_{i-1}<a_{j+1}\delta_j+\delta_{j+1}=\delta_{j-1}$$
  $\blacktriangleright$ Case 2 : if $(-1)^{j+1}(d'_{j+1}-d_{j+1})=a_{j+1}$.\\ 
 $\blacktriangleright\blacktriangleright$  Subcase 1 : if $j$ is even, $d'_{j+1}=0$ and $d_{j+1}=a_{j+1}$.\\
 \esp We can not have $d'_j=a_j$, for $(-1)^{j-1}d_j<(-1)^{j-1}d'_j$, so $d'_i=0$ for all $i>j$ and, since $d$ does not end with $(\max,0)^{\infty}$, then :
   $$ \sum\limits_{i=j+1}^{\infty} (-1)^{i}(d'_i-d_i)\delta_{i-1}=a_{j+1}\delta_j - \sum_{i=j+2}^{\infty}  (-1)^{i}d_i\delta_{i-1}<   \sum_{p=0}^{\infty}  a_{j+2p+1}\delta_{j+2p} $$
 \esp \esp Indeed, $(-1)^{j+2p}d_{j+2p}\delta_{j+2p-1}\geqslant 0$, for all $p\in\N$, since $j$ is even. So :
 $$   \sum\limits_{i=j+1}^{\infty} (-1)^{i}(d'_i-d_i)\delta_{i-1} <   \sum_{p=0}^{\infty} (\delta_{j+2p-1}-\delta_{j+2p+1})= \delta_{j-1}$$
 $\blacktriangleright\blacktriangleright$ Subcase 2 : if $j$ is odd, similar arguments lead to the same conclusion ( we swap $d$ and $d'$).\end{mademo}

$\bullet$ Now, we define two maps on these sets :

 \begin{maprop} .\\
   (i) the map $\Psi_{\alpha}$ ( defined below) is an order isomorphism from $(E_{(\alpha)} ,\leqslant_R)$ to $(\N,\leqslant)$.
     $$ \Psi_{\alpha} : \begin{cases} E_{(\alpha)} \to \N \\ d\to \sum\limits_{j=1}^{\infty}d_jq_{j-1} \end{cases}$$
     
     (ii) the map $\Lambda_{\alpha} $ ( defined below) is an order isomorphism from $(E_{\alpha} ,\leqslant_A)$ to $([0,1[,\leqslant)$. : 
       $$ \Lambda_{\alpha} : \begin{cases} E_{\alpha} \to [0,1) \\ d\to \sum\limits_{j=1}^{\infty}d_j(-1)^{j-1}\delta_{j-1} \end{cases}$$
       
     (iii) we have : 
     $$ \forall n\in \N, \{ n\alpha \}= \Lambda_{\alpha}(\Psi_{\alpha}^{-1}(n))$$  
   \end{maprop}

 \textbf{Remark  1 : } the infinite sum in the definition of $\Psi_{\alpha}$ is in fact a finite one.  The infinite sum in the definition of $\Lambda_{\alpha}$ is well defined since :
 $$ \forall j\in \N^*, 0\leqslant d_{j}\delta_{j-1} \leqslant a_{j}\delta_{j-1} = \delta_{j-2}-\delta_{j}$$

  \textbf{Remark  2 : } if we had defined $E_{\alpha}$ without the restriction about the ending of the sequences, then the result about $\Lambda_{\alpha}$ would have been valid, except that : for $x\in \{  \{ n\alpha \}, n\in\N\}$, $x$ would have  three ( two for $0$) preimages : the one in  $E_{(\alpha)}$ and those that end with $(\max,0)^{\infty}$, an infinite sequence of "$a_k,0$".\\

  \begin{mademo}
  (i) see proof of Lemma 3 and proof of Algorithm 1.\\
  (ii) first, we will prove that $\Lambda_{\alpha}$ is increasing : let $d,d'\in E_{\alpha}$ such that $d<_A d'$. Then, we have $j\in\N^*$ such that :
  $$(-1)^{j-1}d_j<(-1)^{j-1}d'_j \esp \text{ and }  \esp \forall i<j, d_i=d'_i$$
  \esp So :
  $$ \Lambda_{\alpha}(d')-\Lambda_{\alpha}(d)= (-1)^{j-1}(d'_j-d_j)\delta_{j-1} +\sum_{i=j+1}^{\infty} (-1)^{i-1}(d'_i-d_i)\delta_{i-1} $$
  \esp But, $(-1)^{j-1}(d'_j-d_j)\delta_{j-1}\geqslant \delta_{j-1}$, so with Lemma 5, we obtain : $ \Lambda_{\alpha}(d')-\Lambda_{\alpha}(d)>0$.\\
\esp Now that we have proved that $\Lambda_{\alpha}$ is increasing, we can easily deduce that $\Lambda_{\alpha}(E_{\alpha})\subset [0,1[$ : first, remark that $(0,\cdots,0)$ is the lowest element of $E_{\alpha}$ ( with ALO), so $\Lambda_{\alpha}(d)\geqslant 0$ for all $d\in E_{\alpha}$. In addition, if $j$ is even $(-1)^{j-1}d_j\delta_{j-1}\leqslant 0$ and if $j$ is odd, say $j=2p+1$, with $p$ a non negative integer, then $(-1)^{j-1}d_j\delta_{j-1}\leqslant a_{2p+1}\delta_{2p}$, this inequality being strict for at least one $p$, so :

$$ \forall d\in E_{\alpha}, \Lambda_{\alpha}(d)< \sum_{p=0}^{\infty} a_{2p+1}\delta_{2p}=\sum_{p=0}^{\infty}(\delta_{2p-1}-\delta_{2p+1})=\delta_{-1}=1$$
\esp For the surjectivity, we refer to Algorithm 3(ii) below.\\
(iii) see proof of Proposition 3(ii). 
\end{mademo}

\begin{monalgo} .\\
  (i)  the inverse function of $\Psi_{\alpha}$ is defined by the following algorithm : \\
  \esp Let $n\in\N$ and $r=\max(\{ k\in\N, n<q_k+q_{k-1} \})$. We define $d$ by :  $ \forall k>r, d_k =0 $ and \\
%     \begin{algorithm}
%\caption{Voici un algorithme.}
%    \KwIn{$n$}
%    \KwOut{$(d_i)_{i\in \inter{1}{r}}$}
%\BlankLine

%\For{$k \gets r$ \KwTo $1$ \KwStep $-1$}{
% $d_k \gets \max\left(0,\left\lfloor \frac{n-q_{k-2}}{q_{k-1}}%\right\rfloor \right)$ \; 
% $ n \gets n-d_kq_{k-1}$}
%\Return{$x$}
%\end{algorithm}
 \shadowbox{\parbox{11cm}{\textbf{Input : } $n$ \esp \textbf{Output : } $(d_i)_{i\in \inter{1}{r}}$ \\ 
   for $k=r$ to $k=1$ with  step $-1$ : $\begin{cases} d_k=\max\left(0,\left\lfloor \frac{n-q_{k-2}}{q_{k-1}}\right\rfloor \right)\\ n \leftarrow n-d_kq_{k-1} \end{cases}$}}\\
 %  Then $d\in E_{(\alpha)}$ and $\Psi_{\alpha}(d)=n$.\\
 
 (ii) the inverse function  of $\Lambda_{\alpha}$ is defined by the following ( infinite) " algorithm" : \\
  \esp Let $\beta\in [0,1[$. We denote $\beta_0=\beta$ and define the sequences $b=(b_k)_{k\in\N^*}$ and $(\beta_k)_{k\in\N^*}$ by :\\
  
%    \begin{algorithm}
%\caption{Voici un algorithme.}
%    \KwIn{$\beta$}
%    \KwOut{$(b_i)_{i\in \N}$}
%\BlankLine
 
%\For{$k \gets 1$ \KwTo $+\infty$}{
% $b_k \gets \min\left(a_k,\left\lceil \frac{\beta}%{\delta_{k-1}}\right\rceil \right)$ \; 
% $ \beta \gets b_k\delta_{k-1} -\beta$}
%\Return{$x$}
%\end{algorithm}
  \shadowbox{\parbox{11cm}{\textbf{Input : } $\beta$ \esp \textbf{Output : } $(b_i)_{i\in \N^*}$ \\ 
  for $k=1$ to $k=\infty$ with step $1$ : $\begin{cases} b_k=\min\left(a_k,\left\lceil \frac{\beta_{k-1}}{\delta_{k-1}}\right\rceil\right) \\ \beta_{k}=b_k\delta_{k-1}-\beta_{k-1} \end{cases} $}}
 % \esp $b_{r+1}=\beta_rq_{r+1}$.
 % We have then $b\in E_{\alpha}$ where $b=(b_k)_{k\in\N^*}$ and $\beta = \Lambda_{\alpha}(b)$.

 \end{monalgo}
      
\begin{mademo}
(i) see proof of Algorithm 1.\\
(ii) the proof that $b\in E_{\alpha}$ is the same as the proof of Algorithm 2, with the additional argument : $b$ does not end with $(\max,0)^{\infty}$, an infinite sequence of " $(a_k,0)$", that will be shown below. \\
\esp First, we remark that $(\beta_k)_k$ converges towards $0$, for $(-1)^{k}\beta_k-(-1)^{k-1}\beta_{k-1}=(-1)^kb_k\delta_{k-1}$ is the general term of a convergent serie. We can define $\beta'=\sum\limits_{j=1}^{\infty}b_j(-1)^{j-1}\delta_{j-1}$ and verify that $\beta'=\beta$ : 
$$ \beta'= \sum_{j=1}^{\infty}(-1)^{j-1}(\beta_j+\beta_{j-1})=\beta_0=\beta $$
\esp Suppose that $b$ ends with $(\max,0)^{\infty}$ : this means that, we have $r\in \N^*$, such that :
$$ ( r=1 \text{ or } b_{r-1}\not = 0 ) \esp ; \esp  \forall p\in \N, b_{r+2p}=a_{r+2p} \esp ; \esp b_{r+2p+1}=0 $$ 
\esp So :
$$ \beta= \sum\limits_{j=1}^{r-1}b_j(-1)^{j-1}\delta_{j-1}+(-1)^{r-1}\sum_{p=0}^{\infty}a_{r+2p}\delta_{r+2p-1}=\sum\limits_{j=1}^{r-1}b_j(-1)^{j-1}\delta_{j-1}+(-1)^{r-1}\delta_{r-2} $$
\esp If $r=1$, then $\beta=\delta_{-1}=1$, so $r\geqslant 2$ and  we recognize $\beta=\Lambda_{\alpha}(b')$, where $b'=(b_1,\cdots,b_{r-2},b_{r-1}-1)\in E_{(\alpha)}$. Using the proof of Algorithm 2, we obtain $\beta_{r-1}=0$, so $b$ ends with an infinite sequence of " $0$".
   \end{mademo}   
    
    \esp We can sum up these formulae : for all non negative integers $n$, if we denote  $d=\Psi_{\alpha}^{-1}(n)$ :
     $$n=\sum_{j=1}^{\infty} d_jq_{j-1} \esp ;\esp   \lfloor n\alpha\rfloor =\sum_{j=1}^{\infty} d_jp_{j-1}  \esp ; \esp \{ n\alpha \}=\sum_{j=1}^{\infty} (-1)^{j-1}d_j\delta_{j-1}$$

  \textbf{Notations : } if no ambiguity, we will denote $n=(d_1,d_2,\cdots,d_s)_{\alpha}$ the $\Psi_{\alpha}$-numeration of an integer $n$ and $\beta=(b_1,\cdots)_{\alpha}$ the $\Lambda_{\alpha}$-numeration of a real $\beta$ of $[0,1[$.\\

 \textbf{Remark 1 : }  we denote $\N_{\alpha}$ the completion of $(\N,D)$, where $D$ is the distance defined by :
 $$ \forall n,n'\in\N, \esp  D(n,n')= |\{ n'\alpha \}-\{ n\alpha \}|$$
 \esp Proposition 4 proves that $\N_{\alpha}$ can be represented ( bijectively) by $E_{\alpha}$ : if $n \in \N_{\alpha}$ is represented by $d\in E_{\alpha}$ then we could define : $\{ n \alpha\} := \sum\limits_{j=1}^{\infty} (-1)^{j-1}d_j\delta_{j-1}$. We obtain a bijective map :
 $$   \N_{\alpha}\to [0,1[ \esp ; \esp n \to  \{ n \alpha\}$$ 

  \textbf{Remark 2 : }  in next subsection, we will study the effect of the symmetry $\beta\to 1-\beta$ on $\alpha$-numeration of reals of $[0,1[$. But now, we are interested in this symmetry acting both on $\alpha$ and $\beta$, which gives a much simpler result :\\
  - first, let $\alpha$ be a real in $]0,1/2[$ and let us consider the CFE of $\alpha$ and $1-\alpha$ :
  $$ \alpha=[a_k]_{k\in\N} \esp \Rightarrow \esp 1-\alpha=[0,1,a_1-1,a_{[2,\infty]}]$$ 
  \esp Indeed, if we denote $1-\alpha=[a'_k]_{k\in\N},\alpha_1=[a_{[2,\infty]}]$ and $\alpha'_1=[a'_{[2,\infty]}]$, then $a'_0=a_0=0, a'_1=1$ and :
  $$ \alpha=\frac{1}{a_1+\alpha_1} \esp ; \esp 1-\alpha=\frac{1}{1+\alpha'_1}$$
  \esp So :
  $$ \alpha'_1= \frac{1}{1-\alpha}-1= \frac{1}{\frac{1}{\alpha}-1}= \frac{1}{a_1-1+\alpha_1}$$
 - secondly : let $\alpha\in ]0,1/2[,\beta\in ]0,1[$ and $(b_k)_k$ its $\alpha$-numeration, then :
 $$ 1-\beta= (1,b_1-1,b_{[2,\infty]})_{1-\alpha}$$
 \esp Indeed : if we denote $\delta'_i$ the analoguous of $\delta_i$ ( related to $\alpha$) for $1-\alpha$ ( see above), then :
 $$ \delta'_{-1}=1 \esp ; \esp \delta'_0=1-\alpha \esp ; \esp  \forall i\geqslant 1, \esp \delta'_i=\delta_{i-1} $$
 \esp The last equality is obtained with obvious induction and previous result on CFE. Now, we just have to verify that :
 $$ \delta'_0-(b_1-1)\delta'_1+\sum_{i\geqslant 2}(-1)^{i-1}b_{i-1}\delta'_{i-1}=1-\beta$$
 \esp that is an easy calculation...
 \newpage

   \subsection{$\alpha$-numeration of negative integers}
    
     \esp Let $\alpha$ be an irrational in $]0,1[$ and $[a_k]_{k\in\N}$ its CFE.   We have seen at 2.3 that $E_{\alpha}$, the set of $\alpha$-admissible sequences is in bijective correspondance  with $[0,1[$, via the following map :
     $$ \Lambda_{\alpha } \esp : \esp  d=(d_k)_{k\in\N^*} \to \sum_{k=1}^{\infty} d_k\delta'_{k-1}$$
     \esp where $\delta'$ is the sequence defined by :
     $$ \delta'_{-1}=-1 \esp ; \esp \delta'_0=\alpha \esp  ; \esp \forall k\in \N^*, \esp \delta'_k=a_k\delta'_{k-1}+\delta'_{k-2}$$
     \esp with notations of 2.3, we have :
     $$ \forall k\in \inter{-1}{+\infty}, \esp \delta'_k=q_k\alpha-p_k=(-1)^k\delta_k$$
     \esp In addition $\delta'$ converges towards $0$ and we could set $\delta'_{\infty}=0$. \\

     \esp In order to define the $\alpha$-numeration of negative integers, we consider the natural involution of $[0,1[$, that we denote $C$  : the complement to $1$.
     $$ C : \begin{cases} C(0)=0 \\ \forall x\in ]0,1[, C(x)=1-x \end{cases} $$
     \esp We also have : $\forall x\in [0,1[, C(x)=\{ -x \}$.  We can see $C$ as the usual conjugacy over the unit circle $\mathbb{U}$, the set of complex of moduli one,  via the bijection : $[0,1[\to \mathbb{U}, x\to e^{ 2i\pi x}$. $C$ is decreasing, when restricted to $]0,1[$.  \\
     
     \textbf{Question : } is there a simple and natural expression of conjugate involution $C_{\alpha}$ of $E_{\alpha}$, induced by $C$, via $\Lambda_{\alpha }$, that is :
     $$ C_{\alpha}= \Lambda_{\alpha }^{-1}\circ C \circ \Lambda_{\alpha }$$
     \esp Thinking of the analoguous problem for usual $(b^k)_k$ basis-numeration, where $b$ is an integer bigger than $1$, we could try to use a  kind of " complement to $(a_k)_{k\in \N^*}$" transformation. Indeed, $(a_k)_{k\in \N^*}$ is the biggest sequence in $E_{\alpha}$ for the usual lexicographic order. But, we also have to add 1 to the first digit, so, let $m$ be the following sequence :
     $$m_1=a_1+1 \esp ; \esp \forall k>1, \esp m_k=a_k $$
     \esp We extend the definition of $\Psi_{\alpha}$ to all real sequences in $l^1(\delta')=\{ u\in \R^{\N^*}, \sum_k|u_k\delta'_k| < +\infty \}$. 
     $$ L_{\alpha} : l^1(\delta') \to \R ; \esp d \to \sum_{k=1}^{\infty} d_k\delta'_{k-1}$$
     
    \esp Then, $L_{\alpha}(m)=1$, for :
    $$   \sum_{k=1}^{\infty} m_k\delta'_{k-1}=\alpha+ \sum_{k=1}^{\infty} a_k\delta'_{k-1}=\alpha+ \sum_{k=1}^{\infty} (\delta'_k-\delta'_{k-2})=\alpha-\delta'_{-1}-\delta'_0=1 $$
       
      \esp Since $L_{\alpha}$  is linear, we have :
      $$ \forall d\in l^1(\delta'), \esp L_{\alpha}(m-d)=1-L_{\alpha}(d)$$ 
      \esp In particular, for $d\in E_{\alpha}$, we obtain : $L_{\alpha}(m-d)=1-\Psi_{\alpha}(d)$. So, the question is : do we always have $m-d\in E_{\alpha}$ ? Unfortunately, no. But,  $m-d\in E_{\alpha}$ in most cases. \\
      \esp First, since $d$   does not end with $(\max,0)^{\infty}$ ( see 1.2), that is also the case for $m-d$. \\
      \esp Secondly, if $d$ is not the null sequence, then $m_k-d_k\in \inter{0}{a_k}$ for all $k\in \N^*$, and $m_1-d_1>0$.\\
      \esp Finally, the only case where $d\in E_{\alpha}$ and $m-d\not\in E_{\alpha}$  is when  $m-d$ contains a finite word of consecutive $0$, that is not preceeded by a maximal digit ( say $d_k=a_k$) and that is not succeeded by a $0$.  We will name such a word, a \emph{not admissible word}. Such a word can appear in $m-d$, for $d$ can contain a word with consecutive maximal digits.\\
      \esp We will see below how to convert such a sequence into an $\alpha$-admissible sequence. First, let $\sim$ denote the  equivalence relation on $l^1(\delta')$, induced by $L_{\alpha}$ :
      $$ \forall u,v\in l^1(\delta'), \esp u\sim v \Leftrightarrow L_{\alpha}(u)=L_{\alpha}(v)$$
       \esp This relation $\sim$ is compatible with the linear structure of $l^1(\delta')$.\\
   %  \esp Let denote $r$ and $s$ two positive integers  and $B_{r,s}$, the finite word :
   %  $$ B_{r,s}= (a_r,0,a_{r+2},0,\cdots,0,a_{r+2s})$$

   \esp  We have, for all $r,s\in \N^*$  : 
      $$  (0^r,1,(\max,0)^{s-1},\max,-1,0^{\infty}) \sim 0^{\infty}  \esp \esp \mathbf{(1)}$$
     \esp Indeed :
     $$ L_{\alpha}((0^r,1,(\max,0)^{s-1},\max,-1,0^{\infty}))=\delta'_r+\sum_{k=1}^{s}a_{r+2k}\delta'_{r+2k-1}- \delta'_{r+2s}= $$
     $$=\delta'_r+\sum_{k=1}^{s}(\delta'_{r+2k}-\delta'_{r+2k-2})- \delta'_{r+2s}=0$$
     
     \esp \textbf{Case 1 : } a list of an even number of consecutive $0$ ( not preceeded by a maximal digit and not succeeded by a $0$). So, if we have a sequence $(e_k)_k$, such that $e_{[1,r]}=[e_1,\cdots,e_r]$ only contains admissible words and such that $e_r\not = a_r,e_{r+2s+1}\not =0$ and $e_k=0$ for $k\in \inter{r+1}{r+2s}$ ( where $r,s\in\N^*$). \\
      \esp Then, adding $(e_k)_k$ to relation \textbf{(1)}, we obtain  : 
     $$ (e_k)_{k \geqslant 1}\sim (e_{[1,r]},1, (\max,0)^{s-1},\max,e_{r+2s+1}-1,e_{[r+2s+2,\infty]}) $$
      \esp Thus, the new sequence $(e'_k)_k$ only contains admissible words in its first $r+2s+1$ digits.\\
      
        \esp \textbf{Case 2 : } a list of an odd number of consecutive $0$ ( not preceeded by a maximal digit and not succeeded by a $0$). So, if we have a sequence $(e_k)_k$, such that $e_{[1,r]}$ only contains admissible words and such that $e_r\not = a_r,e_{r+2s}\not =0$ and $e_k=0$ for $k\in \inter{r+1}{r+2s-1}$ ( $r,s\in\N^*$). \\
      \esp    \esp Then, adding $(e_k)_k$ to relation \textbf{(1)} ( with $r-1$ instead of $r$), we obtain  : 
     $$ (e_k)_{k\geqslant 1}\sim (e_{[1,r-1]},e_r+1,  (\max,0)^{s-1},\max,e_{r+2s}-1,e_{[r+2s+1,\infty]}) $$
      \esp Thus, the new sequence $(e'_k)_k$ does not contain any not admissible word in its first $r+2s$ digits.\\

      \esp In both cases, we have converted the not admissible word of $(e_k)_k$ into an admissible word, giving the same image for $L_{\alpha}$. This provides a ( possibly infinite) process to convert  any not admissible element of $m-E_{\alpha}$ into an element of $E_{\alpha}$. We only have to browse once the sequence $(e_k)_k$ to convert it into an equivalent $\alpha$-admissible sequence : \\
      
      \textbf{Process of conversion : } \\
    \resultat{  let $d$ denote an $\alpha$-admissible sequence that is not the null sequence and $e=m-d$. Then $e\in \inter{1}{a_1}\times \prod_{k>1}\inter{0}{a_k}$. We denote $(r_j)_j$ and $(s_j)_j$ the sequences of positive integers such that, the finite lists of consecutive $0$ in $e$ are for indices from $r_j+1$ to $r_j+2s_j$ or $r_j+2s_j-1$, depending on the parity of the lengths $(l_j)_j$ of these lists. We apply then the inductive following process : \\
     \esp We suppose that we have converted the digits of $e$ for the indices $k\leqslant r_j$. Then : we can suppose that $e_{r_j}<a_j$ ( if $e_{r_j}=a_j$, then we change $r_j \leftarrow r_j+1$) and $e_{r+l_j+1}>0$.\\
     - Case 1 : if $l_j$ is even, then :
     $$ e \leftarrow (e_{[1,r_j]},1, (\max,0)^{s_j-1},\max,e_{r_j+2s_j+1}-1,e_{[r_j+2s_j+2,\infty]}) $$       
     - Case 2 : if $l_j$ is odd, then :
      $$ e \leftarrow (e_{[1,r_j-1]},e_{r_j}+1, (\max,0)^{s_j-1},\max,e_{r_j+2s_j}-1,e_{[r_j+2s_j+1,\infty]})$$   }\\
      
     \esp So, this process explicits the map $C_{\alpha}$, that is the relation between the $\alpha$-numerations of $\beta$ and $1-\beta$ for a real $\beta\in ]0,1[$. We will name this map : \emph{ CFE-complement}.\\
     \esp Now, let us consider the particular case of $\beta=\{ n\alpha\}$, where $n\in\N^*$. We have seen in 2.3 that $n$ and $\beta$ have the same $\alpha$-numeration. Since $\{ -n\alpha \}=1-\beta$, it is natural to define the $\alpha$-numeration of $-n$ as follows :
     
     \begin{madef} [$\alpha$-numeration of a negative integer].\\
     for any positive integer $n$, we define the $\alpha$-numeration of $-n$ as the CFE-complement of  the $\alpha$-numeration of $n$.
     \end{madef}

      \textbf{Notations : } we denote $E^c_{(\alpha)}$ the subset of $E_{\alpha}$ of sequences ending  with $\max^{\infty}$, that is to say :
      $$ E^c_{(\alpha)}=\{ e\in E_{\alpha},  \exists k\in\N,\forall i>k, \esp e_i=a_i\}$$
      \esp We have then $E^c_{(\alpha)}=C_{\alpha}(E_{(\alpha)})$ and $E^c_{(\alpha)}$ is the set of $\alpha$-admissible sequences that " $\alpha$-numerate" negative integers ( see Proposition below). \\
      \esp We will also denote $F_{\alpha}=E_{(\alpha)}\cup E^c_{(\alpha)}$ and we extend RLO, that we defined on $E_{(\alpha)}$, to $F_{\alpha}$ :
        $$ \forall d,d'\in F_{\alpha}, \esp d<_R d' \Leftrightarrow \exists k\in \N^*, ( (d_k<d'_k, \forall i>k, d_i=d'_i) \text{ or } ( \forall i\geqslant k, d_i=a_i, d'_i=0))$$
      
      \textbf{Remark : } the above process of conversion is, in that frame, an algorithm, since an element of $E^c_{(\alpha)}$ only contains a finite number of lists of consecutive $0$. 
      
      \begin{maprop}  we can extend $\Psi_{\alpha}$ from $E_{(\alpha)}$ to $F_{\alpha}$ as follows :
      $$ \forall e\in E^c_{(\alpha)}, \esp \widetilde{\Psi}_{\alpha}(e)=-1- \sum_{k=1}^{\infty}(a_k-e_k)q_{k-1} $$
      hence, $\widetilde{\Psi}_{\alpha}$ is an order isomorphisme from $(F_{\alpha},\leqslant_R)$ to $(\Z,\leqslant)$ and we still have  : 
      $$ \forall n\in \Z,\esp  \Lambda_{\alpha}(\widetilde{\Psi}_{\alpha}^{-1}(n))=\{ n\alpha \}$$

    \end{maprop}

    \begin{mademo}
    - Formula and injectivity : let $e\in E^c_{(\alpha)}$. First, we remark that the sum in the definition of $\widetilde{\Psi}_{\alpha}(e)$ is finite, since $e_k=a_k$ for $k$ large enough. Let denote $d=m-e$ and :
    $$ n=\sum_{k=1}^{\infty} d_kq_{k-1} \esp ; \esp   \beta=\sum_{k=1}^{\infty}d_k\delta'_{k-1}$$
       \esp Now, $ \widetilde{\Psi}_{\alpha}(e)=-n$ and $n$ is a positive integer ( $d$ ends with $0^{\infty}$ and $d_1>0$), so $ \widetilde{\Psi}_{\alpha}(e)\in \Z_-^*$.\\
       \esp  But $d$ is not always in $E_{(\alpha)}$. Nevertheless $\beta\in ]0,1[$ for :
    $$\beta=L_{\alpha}(d)=L_{\alpha}(m)-L_{\alpha}(e)=1-\Lambda_{\alpha}(e)$$
    \esp indeed, $e\in E_{\alpha}$. Finally :
    $$ n\alpha-\beta=\sum_{k=1}^{\infty}d_kp_{k-1} \in \N$$
    \esp so : $\beta=\{ n\alpha \}$. We obtain : $\Lambda_{\alpha}(e)=1-\beta=\{ -n\alpha\}$. We can conclude :
    $$ \forall e\in F_{(\alpha)}, \esp \Lambda_{\alpha}(e)=\{ \widetilde{\Psi}_{\alpha}(e) \alpha \}  \esp (1)$$
    \esp Since, $\Lambda_{\alpha}$ is injective, we deduce that $\widetilde{\Psi}_{\alpha}$ is injective. \\
    - Surjectivity : let $n\in \N^*,d=\Psi_{\alpha}^{-1}(n)$ and $e=\Lambda_{\alpha}^{-1}(1-\{ n\alpha \})$. Then, $e\in E^c_{(\alpha)}$ ( see the beginning of this section) and, with (1) :
    $$ \Lambda_{\alpha}(e)=\{ \widetilde{\Psi}_{\alpha}(e) \alpha \}  =1-\{ n\alpha \}=\{ -n\alpha \}$$
    \esp So, $ \widetilde{\Psi}_{\alpha}(e)=-n$, for $\widetilde{\Psi}_{\alpha}(e)\in \Z$. So, $\widetilde{\Psi}_{\alpha}$ is surjective. \\
    - Increase : let $e,e'\in E^c_{(\alpha)}$ such that $e<_R e'$. \\
    --- Case 1 : if $e\in E^c_{(\alpha)}$ and $e'\in E_{(\alpha)}$, then $\widetilde{\Psi}_{\alpha}(e)<0\leqslant \widetilde{\Psi}_{\alpha}(e')$. \\
    --- Case 2 : if $e,e'\in E_{(\alpha)}$, we have proved in Proposition 2 that $\Psi_{\alpha}(e)< \Psi_{\alpha}(e')$. \\
    --- Case 3 : if $e,e'\in E^c_{(\alpha)}$, then :
    $$ \widetilde{\Psi}_{\alpha}(e')-\widetilde{\Psi}_{\alpha}(e)=\sum_{k=1}^{\infty}(e'_k-e_k)q_{k-1}= \Psi_{\alpha}(d')-\Psi_{\alpha}(d)$$
    \esp where $d=((e_k)_{k\in \inter{1}{r}},0^{\infty})$ and $d'=((e'_k)_{k\in \inter{1}{r}},0^{\infty})$, the integer $r$ being such that $e'_i=e_i=a_i$ for $i>r$. Since $e,e'$ are $\alpha$-admissible, we can claim that $d,d'\in E_{(\alpha)}$. So, with Proposition 2, $\Psi_{\alpha}(d')-\Psi_{\alpha}(d)>0$. So, $\widetilde{\Psi}_{\alpha}(e')-\widetilde{\Psi}_{\alpha}(e)>0$.\\
    \esp We have proved that $\widetilde{\Psi}_{\alpha}$ is increasing on $F_{\alpha}$.
       \end{mademo}

        \esp Note that the definition of $\widetilde{\Psi}_{\alpha}$ in Proposition 5 could be given by the same formula for $d$ in $E_{(\alpha)}$ and for $d$ in $ E^c_{(\alpha)}$, with the following convention : $+\infty=0$, so that $q_{n}\tend 0$. Indeed, if we define :
        $$ \forall d\in F_{\alpha}, \esp \Psi_{\alpha}(d)=\sum_{k=1}^{\infty}d_kq_{k-1}$$
        \esp then, it is convenient, since :
        $$\sum_{k=1}^{\infty}a_kq_{k-1}=\sum_{k=1}^{\infty}(q_k-q_{k-2})=0+0-q_0-q_{-1}=-1$$
       \esp We also have, with this convention a coherent result for both " improper expansions" of an integer $n$, herited from improper expansions of $\{ n\alpha \}$ ( see remark 2, below Proposition 4), whose proper expansion is $(d_1,d_2,\cdots,d_r)$ with $d_r>0$. Indeed, these improper expansions are $(d_{[1,r]},1,(\max,0)^{\infty})$ and $(d_{[1,r-1]},d_r+1,(\max,0)^{\infty})$ ( if $d_r<a_r$) or $(d_{[1,r]},0,1,(\max,0)^{\infty})$ ( if $d_r=a_r$). Moreover :
       $$ \forall s\in \N, \esp \sum_{j=0}^{\infty} a_{s+2j+1}q_{s+2j}= \sum_{j=0}^{\infty}(q_{s+2j+1}-q_{s+2j-1})=0-q_{s-1}$$

   \section{Complements}
   \label{ch:dyn}

 \subsection{dynamic  generating $\alpha$-numeration}

$\bullet$ What follows is inspired by the analoguous result for the usual Ostrowski numeration made by Ito in  \textbf{[5]} :

  \begin{maprop} let $\alpha$ be an irrational and $[a_k]_{k\in\N}$ its CFE. Let $\beta\in [0,1[$ and $(b_k)_k$ its $\alpha$-numeration. We have : 
  $$  \forall k\in \N^*,\esp (a_k,b_k)=AH^{k-1}(\alpha,\beta)$$
  \esp where $H$ is a self map of the open trapezoid $U$ defined by : for $(x,y)\in\R^2$
  $$(x,y)\in U \Leftrightarrow \begin{cases}  0<x<1 \\ -x<y<1  \end{cases} $$
  $$  H : (x,y)\to \left( \left\{\frac{1}{x}\right\},\min\left( \left\lfloor \frac{1}{x} \right\rfloor, \left\lceil \frac{y}{x} \right \rceil \right) - \frac{y}{x} \right)$$
   $$  A : (x,y)\to \left( \left\lfloor \frac{1}{x}\right\rfloor,\min\left( \left\lfloor \frac{1}{x} \right\rfloor, \left\lceil \frac{y}{x} \right \rceil \right)  \right)$$

  \end{maprop}
  
  \textbf{Remark 1 : } we could prefer the following expressions, distinguishing two cases :
  $$ \forall (x,y)\in U, \esp \begin{cases}  A(x,y)= (\lfloor 1/x \rfloor, \lceil y/x \rceil ) \esp ; \esp  H(x,y)=(\{ 1/x\}, \{ -y/x \})\text{ if } y\leqslant x\lfloor 1/x \rfloor \\
  A(x,y)= (\lfloor 1/x \rfloor, \lfloor 1/x \rfloor ) \esp ; \esp H(x,y)=(\{ 1/x\}, \{ -y/x \}-1) \text{ else }\end{cases}$$ 
  \esp Indeed, if $y>x\lfloor 1/x \rfloor$, then : $\lfloor 1/x \rfloor < y/x< 1/x$, so $\lfloor 1/x \rfloor=\lceil y/x \rceil-1$.\\
  
  \textbf{Remark 2 : } let us verify that $H(U)\subset U$ : if $y\leqslant x\lfloor 1/x \rfloor $, that is obvious. Else, $ \{ -y/x \}-1= -\{ y/x \}>-\{ 1/x\}$, for $\lfloor 1/x \rfloor < y/x< 1/x$ and so $\{ y/x \}< \{ 1/x\}$ ( see remark 1).\\

 \begin{mademo} we denote $(\alpha_k,\gamma_k)=H^k(\alpha,\beta)$ for all $k\in \N$. We avoid here the notation $\beta_k$  for it is used below as reference to Algorithm 3.\\
 \esp   We already know that $a_k=p_x(AH^{k-1}(\alpha,\beta))$, where $p_x:(x,y)\to x$, since $T(x)=p_x(H(x,y))$ for all $x,y\in ]0,1[$ ( $T$ is the Gauss map, see 1.3). By definition, we have :
 $$ \gamma_0=\beta \esp ; \esp \forall k\in\N^*, \esp \gamma_k=\min(a_k,\lceil \gamma_{k-1}/\alpha_{k-1} \rceil ) -\frac{\gamma_{k-1}}{\alpha_{k-1}}$$
 \esp We denote $\gamma'_k=\frac{\beta_k}{\delta_{k-1}}$, with notations of Algorithm 3 ( see 2.3). We also have :
 $$ \forall i\in\N, \esp \alpha_i=T^i(\alpha)=\frac{\delta_i}{\delta_{i-1}} \esp \text{ so } \esp \frac{\gamma'_i}{\alpha_i}=\frac{\beta_i}{\delta_i}$$
 \esp Thus, according to Algorithm 3 on reals :
 $$ \forall k\in\N^*, \esp b_k=\min(a_k,\lceil \beta_{k-1}/\delta_{k-1} \rceil ) \esp ; \esp \beta_k=b_k\delta_{k-1}-\beta_{k-1}$$
 \esp We deduce :
 $$  \forall k\in\N^*,\esp \gamma'_k= b_k-\frac{\gamma'_{k-1}}{\alpha_{k-1}}$$
 \esp Yet, $\gamma'_0=\beta=\gamma_0$ and we obtain, by obvious induction : $\gamma_k=\gamma'_k$ for all integer $k\in\N$. Then :
 $$ \forall k\in\N^*, \esp b_k=\min(a_k,\lceil \gamma_{k-1}/\alpha_{k-1} \rceil ) $$
 \esp This ends the proof.
 \end{mademo} 
      
 \newpage

   \subsection{$\alpha$-germs and orbits of $\alpha$-rotation}
   
   \esp Our $\alpha$-numeration is related to $f_{\alpha}$, the rotation on the circle $\R/\Z$ defined by :
   $$ \forall x\in \R/\Z,  \esp f_{\alpha}(x)=\alpha+x \esp (\text{ mod } 1)$$
   \esp Let $\alpha$ be an irrational and $[a_k]_k$ its CFE. We know that $f_{\alpha}$ is topologically transitive : its orbits are dense in $X=\R/\Z$. Moreover, it is uniquely ergodic : there is only one   $f_{\alpha}$-invariant ( and ergodic) measure on $X$ : the Lebesgue measure. \\
   \esp Now, we will explicit the conjugate of $f_{\alpha}$ on $E_{\alpha}$, namely the map $g_{\alpha}: E_{\alpha}\to E_{\alpha}$, such that : 
   $$\Lambda_{\alpha}\circ g_{\alpha}=f_{\alpha}\circ \Lambda_{\alpha}$$
   \esp  We remind some notations : $E_{\alpha}$ is the set of $\alpha$-admissible sequences and $F_{\alpha}= E_{(\alpha)}\cup  E^c_{(\alpha)}$, where
   $$ E_{(\alpha)}=\{ (d_k)_k\in E_{\alpha}, \exists n\in \N, \forall k>n,d_k=0 \} \esp ; \esp E^c_{(\alpha)}=\{ (d_k)_k\in E_{\alpha}, \exists n\in \N, \forall k>n,d_k=a_k \}$$
  \esp We will  use an equivalence relation on $E_{\alpha}$, that  defines the notion of \emph{germ} of a sequence :
    $$ \forall d,d'\in E_{\alpha}, \esp d\eqsim d' \Leftrightarrow \exists k\in\N, \forall i>k, d_i=d'_i$$
      \esp We remark that the class of $(0)$ is $E_{(\alpha)}$ and that the class of $(a_k)_{k\in\N^*}$ is $E^c_{(\alpha)}$.\\
  \esp More generally, we can extend RLO to each class of germs of $E_{\alpha}$, as follows :
    $$ (d_k)_k<_R (d'_k)_k \Leftrightarrow \exists j\in \N^*, \begin{cases}  d_j<d'_j \\ \forall k>j, d_k=d'_k \end{cases} $$
  
  \textbf{Remark : } for each class of germs of $E_{\alpha}$, RLO is a total order and every element of the class has a successor ( except for $E^c_{(\alpha)}$, where $(a_k)_{k\in\N^*}$ is the maximal element) and a predecessor ( except for $E_{(\alpha)}$, where $(0)$ is the minimal element).\\
  
 $\bullet$ The following Proposition explicits the orbits of $g_{\alpha}$. Before that, we remark that : for $\beta,\beta'\in \R/\Z$, $\beta$ and $\beta'$ are in the same orbit of $f_{\alpha}$ if and only if it exists $n\in\Z$, such that $\beta'-\beta=n\alpha$ mod $1$. So, an orbit of $g_{\alpha}$ is the set of  $\alpha$-numerations of the $\{ \beta+n\alpha \}, n\in \Z$, for some $\beta\in [0,1[$.
   
   \begin{maprop}Let $\alpha$ be an irrational, $[a_k]_k$ its CFE and $g_{\alpha}$ defined as above, then : \\
   (i) the orbits of $g_{\alpha}$ are exactly the classes of germs of $E_{\alpha}$, except for  the orbit of $(0)$, that is $ F_{\alpha}$.\\
 (ii)  $g_{\alpha}$ is the successor map on each of theses classes ( with RLO).
   
   \end{maprop}
   \begin{mademo}
   \esp First, the class of $(0)$, via $g_{\alpha}$, is $F_{\alpha}$, the set of   $\alpha$-numerations of the $\{ n\alpha \}, n\in \Z$, as we have seen in previous subsection 3.1.\\
   \esp Let $\beta\in [0,1[$ such that $\beta\not\in\{ \{ n\alpha \}, n\in \Z\}$. We denote $b=(b_k)_k$ its $\alpha$-numeration and $C$ the class of germ of $b$ in $E_{\alpha}$. \\
   \esp If $b'\in C$, then we have an integer $r\in\N$, such that $b'_i=b_i$ for all integer $i>r$. We denote $\beta'= \Lambda_{\alpha}(b')$, then :
   $$ \beta'-\beta= \sum_{k=1}^r (b'_k-b_k)\delta'_{k-1}$$
   \esp but, $\delta'_i=\alpha q_i-p_i$ and $q_i,p_i$ are integer for all  $i\in\N$. So, $\beta'-\beta\in \Z+\alpha\Z$ and we conclude that $\beta'$ is in the $f_{\alpha}$-orbit of $\beta$ and that  $b'$ is in the $g_{\alpha}$-orbit of $b$. \\
     \esp Conversely, suppose that $b'$ is in the $g_{\alpha}$-orbit of $b$. We want to show that $b$ and $b'$ have the same germ. By obvious induction, it suffices to show that this is the case for $b'=g_{\alpha}(b)$, that is to say for $\beta'=\beta+\alpha$. But, since $b$ is not $(a_k)_k$, then there exists an index $r$ such that $b_r<a_r$. We denote $d=(b_{[1,r]},0^{\infty})$. Then, the successor of $d$ in $(E_{(\alpha)},RLO)$ is $d'$ such that $d'_i=0$ for all $i>r$. We claim now that $b'=(d'_{(1,r]},b_{[r+1,\infty]})$. Indeed, $b'\in E_{\alpha}$ and :
     $$ \Lambda_{\alpha}(b')-\Lambda_{\alpha}(b)=\Lambda_{\alpha}(d')-\Lambda_{\alpha}(d)=\alpha$$
     \esp So, $b'$ and $b$ have the same germ.\\
     \esp By the way, we have also proved that $g_{\alpha}$ is the successor map on the class of germ of $b$.   
   \end{mademo}
   
   \textbf{Remark 1 : } this proves that $\R/(\Z+\alpha \Z)$ is represented, via our $\alpha$-numeration $\Lambda_{\alpha}$, by germs of sequences of $E_{\alpha}$.\\  
   
   \textbf{Remark 2 : }  we can define, on each orbit $X$ of $f_{\alpha}$, a natural order, which makes them isomorphic to $(\Z,\leqslant)$ ( but not canonically) :
   $$ \forall x,x' \in X, x\leqslant x' \Leftrightarrow \exists n\in \N, x'=f_{\alpha}^n(x) $$
   \esp In the same way, each class of germ of $(E_{\alpha},RLO)$ ( except for the class of $(0)$, where we consider $F_{\alpha}$) is isomorphic to $(\Z,\leqslant)$.\\

%$\bullet$ We define the usual notion of \emph{( diophantine) equivalence} of two reals $\alpha$ and $\alpha'$ : they are equivalent if and only if their partial quotients are equal except for a finite number of them : 
%$$ \alpha \sim \alpha' \Leftrightarrow \exists r\in \N,  T^r(\alpha)=T^r(\alpha') $$ 
%\esp where $T$ is the Gauss map. We also have :
%$$    \alpha \sim \alpha' \Leftrightarrow \exists M\in SL_2(\Z), \esp (\alpha',1) \simeq M(\alpha,1) \text{ in } P(\R^2)$$
%\esp That is to say :
%$$   \alpha \sim \alpha' \Leftrightarrow \exists a,b,c,d \in \Z, \esp ad-bc=1, \alpha'=\frac{a\alpha+b}{c\alpha+d}$$
% \esp Geometrically, the class of equivalence of $\alpha$ is the set of inverse slopes of  the lines $(OM)$, where $M$ is a point of the lattice generated by $(1,0)$ and $(\alpha,1)$.\\
% \esp Thus, $E_{\alpha}$ and $E_{\alpha'}$ have the same germs of sequences if and only if $\alpha$ and $\alpha'$ are equivalent. No surprises with that fact, since $\Z+\alpha\Z$ and $\Z+\alpha\Z'$ are isomorphic if $\alpha$  and $\alpha'$ are equivalent : these are projections on $x$-axis of the lattices $L=(1,0)\Z+ (\alpha,1)\Z$ and $L'=(1,0)\Z+ (\alpha',1)\Z$.
 
 %??????????????????????????????

   $\bullet$ Now, we define, for any $x$ in $\R, ||x||$, the distance of $x$ to $\Z$. We also have : $||x||=\min(\{ x \},\{ -x \})$.  Later, we define several maps on $\R$ by :  for all $\beta\in \R$
   $$ D_{\alpha}( \beta )= \underset{n\to + \infty}{\liminf}(n||n\alpha - \beta ||) \esp ; \esp  D^+_{\alpha}( \beta )= \underset{n\to + \infty}{\liminf}(n\{ n\alpha - \beta \}) \esp ; \esp  D^-_{\alpha}( \beta )= \underset{n\to + \infty}{\liminf}(n \{ \beta-n\alpha \})$$

 %  \esp  We will explore and generalize these maps later in 5.2., but we can make some  remarks : \\ 
    
  \textbf{Remark 3 : } $D_{\alpha}=\min(D^+_{\alpha},D^-_{\alpha})$, for $\liminf$ " respects" the $\min$. \\
 
    \textbf{Remark 4 : } these 3 maps are $f_{\alpha}$-invariant. Indeed, if $x\in \R$, then :
    $$ \forall n\in\N^*,\esp n\{ n\alpha-(x+\alpha) \}=\{ (n-1)\alpha - x \}=\frac{j+1}{j}\times j \{ j\alpha-x \}$$
    \esp where $j=n-1$. But, $\frac{j+1}{j}$ converges to 1 as $j$ tends to infinity, so the $\liminf$ is the same... \\
    \esp This proves that these maps could be defined on $\R/(\Z+\alpha\Z)$, the additive group of orbits of $f_{\alpha}$ and so they only depend on the germ of the $\alpha$-numeration of $\beta\in \R/\Z$. In other words, these maps only depand on the asymptotic behaviour of the $\alpha$-numeration of $\beta$. \\

 % \esp We can ask many questions about these maps : are they constant, finite, null, bounded, completely discontinuous, Lebesgue-measurable, almost everywhere null or infinite. What is their image ? \\
  \esp It is well known that $D_{\alpha}(0)$ is null if and only if the sequence of partial quotients of $\alpha$ is unbounded and that $D_{\alpha}(0)$  can be defined, restricting $n$ to the denominators of convergents of $\alpha$. But, we have more precise results : 
  $$ \underset{n\to + \infty}{\liminf}\left(\frac{1}{a_n+2}\right) \leqslant D_{\alpha}(0) \leqslant \underset{n\to + \infty}{\liminf}\left(\frac{1}{a_n}\right)$$
  \esp Moreover, Dirichlet's theorem on diophantine approximation gives ( see \textbf{[4]}) :
  $$ \forall \beta\in \R, \esp D_{\alpha}(\beta)\leqslant 1  $$
  \esp And Minkowski has proved that ( see \textbf{[4]} again) :
  $$ \forall \beta\in \R\backslash (\Z+\alpha\Z), \esp \min( D_{\alpha}(\beta),D_{\alpha}(1-\beta))\leqslant \frac{1}{4}$$
  \esp In 4.3, we give some results that helps to compute $D_{\alpha}^+(\beta)$ and $D_{\alpha}^-(\beta)$, in relation to the $\alpha$-numeration of $\beta$.

   \newpage

  \subsection{shift and  inductive structure}

  $\bullet$ Let $\alpha$ be a real in $[0,1[$ and $[0,a_1,a_2,\cdots]$ its CFE. We denote $a=(a_1,\cdots)$ and $\sigma$ the usual shift on sequences. We have seen that : if $\alpha$ is not null, then
 $[0,\sigma(a)]$ is the CFE of $T_1(\alpha)$, where $T_1$ is an extension of the Gauss map, described in 1.3.  We recall that $\mu(\alpha)=+\infty$ if $\alpha$ is not rational and $\mu(\alpha)=r$ if $\alpha$ is rational and its CFE is $[0,a_1,\cdots,a_r,1]$. We define inductively the sequence : $(\alpha_k)_k$ by : 

 $$ \alpha_0=\alpha \esp ; \esp \forall k\in \inter{1}{\mu(\alpha)},\esp \alpha_k=\left\{ \frac{1}{\alpha_{k-1}}\right\}$$

\esp With the remark above, we obtain :
 $$ \forall k\in\inter{0}{\mu(\alpha)-1}, \esp \alpha_k=[0,a_{k+1},\cdots] =[0,\sigma^k(a)]$$
  \esp Moreover, if $\alpha$ is rational and $r=\mu(\alpha)$, then $\alpha_r=0$, for $\alpha_{r-1}=[0,a_r,1]=\frac{1}{a_r+1}$. \\
 
  \esp According to the definition of the sets $(E_{\alpha_k})_k$, we can claim :
  $$ \forall b \in E_{\alpha}, \forall k\in\N, ( \sigma^k(b)\in E_{\alpha_k} \Leftrightarrow b_{k+1}\not = 0 \text{ or } \sigma^k(b)=(0) ) $$
  
   \esp In particular :
   $$E_{T(\alpha)}\subset \sigma(E_{\alpha}) \esp \text{ and } \esp \sigma(E_{\alpha}) \backslash E_{T(\alpha)}= \{ 0 \} \times (E_{T^2(\alpha)}\backslash \{(0)\})$$
   \esp  In addition, if we denote for any $k\in \inter{0}{a_1-1}$ :\\
    -  $E_{\alpha,k}$ : the set of $\alpha$-admissible sequences whose first digit is $k$. We have $E_{\alpha,0}=\{ (0)\}$.\\
    -  $E_{\alpha,a_1}$  : the set of $\alpha$-admissible sequences whose first digit is $a_1$ and second is non null, except for $(a_1,0,0,\cdots)$, that is in this set.\\
    -  $E'_{\alpha,a_1}$ :  the set of $\alpha$-admissible sequences whose first digit is $a_1$ and second is null, except for $(a_1,0,0,\cdots)$, that is not in this set.\\
    
    \esp $(E_{\alpha,k})_{k\in \inter{0}{a_1}}\cup E'_{\alpha,a_1}$ is clearly a partition of $E_{\alpha}$ and ALO induces an order on these subsets : ( where $B<_A B'$ means that for every $b\in B$ and $b'\in B'$, we have $b<_A b'$) 
    $$ E_{\alpha,0} <_A E_{\alpha,1} <_A E_{\alpha,2} <_A \cdots <_A E_{\alpha,a_1} <_A E'_{\alpha,a_1}$$
    
    \begin{monlem}.\\
    (i) for any $k\in \inter{1}{a_1}$, the map ( see below) 
     is a bijective decreasing map ( induced by $\sigma$).
     $$ \sigma_k : \begin{cases}(E_{\alpha,k},\leqslant_A) \to (E_{T(\alpha)},\leqslant_A) \\ (k,d_{[2,\infty]})\to (d_{[2,\infty]})\end{cases}$$
    (ii) the map ( see below) is a bijective increasing map ( induced by $\sigma^2$).
     $$ \sigma^{(2)} : \begin{cases}(E'_{\alpha,a_1},\leqslant_A) \to (E_{T^2(\alpha)}\backslash \{ (0) \},\leqslant_A) \\ (a_1,0,d_{[3,\infty]})\to (d_{[3,\infty]})\end{cases}$$
    \end{monlem}
    \begin{mademo}direct consequence of former remarks and definition of sets $E_{\alpha}$ and ALO.    
    \end{mademo}
   
      \esp So to say, $(E_{\alpha},<_A)$ consists in one null element, followed by $a_1$ ordered copies of $(E_{T(\alpha)},<_{A'})$ and, at the end a copy of $(E_{T^2(\alpha)}\backslash \{ (0) \},<_A)$, where $<_{A'}$ denotes inversed ALO.\\

      \newpage
      
      \esp We deduce a result on Kronecker sequences :
      
      \begin{moncoro}
         let $\alpha$ be a real in $[0,1[$, $T$ the usual Gauss map $x\to \{ 1/x \}$.\\ We denote $a_1=\lfloor 1/\alpha \rfloor $ and $K_{\alpha}=\{ \{ k\alpha \}, k\in \N \}$. \\
         \esp The following union are disjoint : 
         $$ K_{\alpha}= \alpha \left( \{ 0 \} \cup \bigcup_{ j\in \inter{1}{a_1 } }( j- K_{T(\alpha)}) \cup (a_1+T(\alpha)(K_{T^2(\alpha)}\backslash \{ 0 \})\right)$$
       \end{moncoro}  
      \begin{mademo} direct consequence of Lemma 6 \end{mademo}

  $\bullet$ Now, we would like to specify the effect of the shift on the integers and reals of $[0,1[$, via their $\alpha$ or $T(\alpha)$-numerations.\\
  
  \esp We  define a sequence of integers $(\nu_k)_k$ by :   
 $$ \nu_0=\nu \esp ; \esp \forall k\in \inter{1}{\mu(\alpha)-2}, \esp \nu_k= \begin{cases} \lfloor \nu_{k-1}\alpha_{k-1}\rfloor  \text{ if } n_{k+1}\not = 0 \text{ or } \sigma^k(n)=(0)\\ 
   \lfloor \nu_{k-1}\alpha_{k-1}\rfloor +1 \text{ else } \end{cases}  $$

 \vspace*{0.5cm}
 
  \begin{monlem} let $k\in\inter{0}{\mu(\alpha)-2}$ and $n=(n_i)_i$ the $\alpha$-numeration of $\nu$ ( we denote $\nu=(n)_{\alpha}$ for example these numeration...) \\
     $\triangleright$ Case 1 : if $n_{k+1}\not = 0$ or $\sigma^k(n)=(0)$, then $\nu_k=(n_{[k+1,\infty]})_{\alpha_k}=\sigma^k(n)_{\alpha_k}$. \\
   $\triangleright$ Case 2 : else $\nu_k=(1,n_{[k+2,\infty]})_{\alpha_k}=(1,\sigma^{k+1}(n))_{\alpha_k}$. 
   \end{monlem}
 
   \begin{mademo}
   we will denote $p_k(x)$ and $q_k(x)$ for the  reduced of the $k^{th}$ convergent of a real $x$, for any non negative integer $k$ and $[a_0(x),a_1(x),\cdots,a_k(x),...]$ its CFE. We have remarked that, if we denote $T(x)=\left\{ \frac{1}{\{ x \}}\right\}$, then :
   $$ \forall j\in\N^*, a_j(T(x))=a_{j+1}(x)$$
   \esp By obvious induction, we can deduce 
    that :
   $$ \forall x\in [0,1[, \forall j\in\N, \esp q_{j-1}(T(x))=p_{j}(x ) \esp \esp (1)$$
   \esp We denote $r=\mu(\alpha)$. Now, we will use an induction on $k\in\inter{0}{r-2}$. Result (i) is true for $k=0$ ( we are in Case 1) . Suppose it is true for $k-1$, where $k\in\inter{1}{r-2}$, then :
   $$ \nu_{k-1}=(n'_k,n_{[k+1,r]})_{\alpha_{k-1}}=(n'_k,\sigma^k(n))$$
   \esp  with $n'_k=1$ or $n_k$, but $n'_k>0$ in all cases. \\
   $\blacktriangleright$ Case 1 : if $n_{k+1}\not = 0$ or $\sigma^k(n)=(0)$, then : with the formula that follows the proof of Algorithm 3 and $(1)$ :
   $$ \nu_k=\lfloor \nu_{k-1}\alpha_{k-1}\rfloor=n'_kp_0(\alpha_{k-1})+\sum_{j=k+1}^rn_jp_{j-k}(\alpha_{k-1})=\sum_{j=k+1}^rn_jq_{j-k-1}(\alpha_{k})$$
   \esp For $p_0(\alpha_{k-1})=0$. So we obtain the $\alpha_k$-numeration of $\nu_k$ : it is $\sigma^k(n)$ for $\sigma^k(n) \in E_{(\alpha_k)}$.     
 %  $$ \nu_k= \sum_{j=k+1}^rn_jq_{j-k-1}(\alpha_{k})$$
  
   $\blacktriangleright$ Case 2 : if $n_{k+1}=0$ and $n_{k+2}\not = 0$, then $(n_{[k+1,\infty]})\not\in E_{\alpha_k}$, but :
   $$ \nu_k=q_0+ \sum_{j=k+2}^rn_jq_{j-k-1}(\alpha_{k})$$
   \esp So  we obtain the $\alpha_k$-numeration of $\nu_k$ : it is $(1,n_{[k+2,\infty]})$ for it is in $E_{(\alpha_k)}$.
 \end{mademo}

 \newpage
 
  \esp We also define a sequence $(\gamma_k)_k$ of reals : 
 $$ \gamma_0=\beta \esp ; \esp \forall k\in \inter{1}{\mu(\alpha)}, \gamma_k= \frac{1}{\alpha_{k-1}}(b_k\alpha_{k-1}-\gamma_{k-1})$$
 
  \vspace*{0.5cm}
  
    \begin{monlem} let $k\in\inter{0}{\mu(\alpha)-2}$. \\
  $\triangleright$ Case 1 : if $b_{k+1}\not = 0$ or $\sigma^k(b)=(0)$, then $\gamma_k=(b_{[k+1,\infty]})_{\alpha_k}=\sigma^k(b)_{\alpha_k}$. \\
   $\triangleright$ Case 2 : else $\gamma_k<0$ and $\gamma_{k+1}=(b_{[k+2,\infty]})_{\alpha_{k+1}}=\sigma^{k+1}(b)_{\alpha_{k+1}}$. 
    
   \end{monlem}
   \begin{mademo} we will use same notations as in previous proof. First, we  remark that ( by obvious induction) :
   $$ \forall x\in ]0,1[, \forall i\in \N, \esp q_i(x)=a_1(x)q_{i-1}(T(x))+p_{i-1}(T(x)) \esp  ; \esp p_i(x)=q_{i-1}(T(x))$$

    \esp We denote $r=\mu(\alpha)$ and argue with induction on $k$. It is clear for $k=0$. Suppose it is true for $k-1$, with $k\in \inter{1}{r-2}$. \\
   --- if $b_k\not = 0$ or $\sigma^{k-1}(b)=(0)$, then $\gamma_{k-1}=(b_{[k,\infty]})_{\alpha_{k-1}}$. So :
   $$ \gamma_{k-1}=\sum_{j=k}^rb_j[\alpha_{k-1} q_{j-k}(\alpha_{k-1})-p_{j-k}(\alpha_{k-1})]=\sum_{j=k}^rb_j[\alpha_{k-1} (a_kq_{j-k-1}(\alpha_{k})+p_{j-k-1}(\alpha_{k}))-q_{j-k-1}(\alpha_k)]$$
   \esp The term of the above sum for $j=k$ is equal to $b_k\alpha_{k-1}$, so :
   $$ \gamma_k= \frac{1}{\alpha_{k-1}}(b_k\alpha_{k-1}-\gamma_{k-1})= \sum_{j=k+1}^rb_j\left[\frac{q_{j-k-1}(\alpha_k)}{\alpha_{k-1}}- (a_kq_{j-k-1}(\alpha_{k})+p_{j-k-1}(\alpha_{k}))\right]$$
   \esp But, $\frac{1}{\alpha_{k-1}}=a_k+\alpha_k$, so :
   $$ \gamma_k= \sum_{j=k+1}^rb_j [\alpha_kq_{j-k-1}(\alpha_{k})-p_{j-k-1}(\alpha_{k})]$$
   \esp Case 1 : $b_{k+1}\not = 0$ or $b_{k+2}=0$ : we recognize the $\alpha_k$-numeration of $\gamma_k$, since $(b_{[k+1,r]})\in E_{\alpha_k}$, with our hypothesis. \\
   \esp Case 2 : $b_{k+1}=0$ and $b_{k+2}\not = 0$, then $(b_{[k+1,r]})\not\in E_{\alpha_k}$ and :
    $$ \gamma_k= \sum_{j=k+2}^rb_j [\alpha_kq_{j-k-1}(\alpha_{k})-p_{j-k-1}(\alpha_{k})]$$
    \esp so :
    $$ \gamma_{k+1}= -\frac{\gamma_k}{\alpha_k}= 
    \sum_{j=k+2}^rb_j \left[ \frac{ p_{j-k-1}(\alpha_{k})}{\alpha_k}-q_{j-k-1}(\alpha_{k})\right]=\sum_{j=k+2}^rb_j[\alpha_{k+1}q_{j-k-2}(\alpha_{k+1})-p_{j-k-2}(\alpha_{k+1})]$$
    \esp the last equality is obtained as above in Case 1...\\
    \esp  Now, $(b_{[k+2,r]}) \in E_{\alpha_{k+1}}$ and $\gamma_{k+1}=(b_{[k+2,r]})_{\alpha_{k+1}}$. We deduce that $\gamma_{k+1}\in ]0,1[$ and $\gamma_k<0$. \\

 --- if $b_k=0$ and $b_{k+1}\not = 0$, then, with induction hypothesis, we obtain the result since we are in Case 1.   
 \end{mademo}

\newpage
      
 \section{Order properties of Kronecker sequences}
\label{ch:kro} 
%  \esp We will denote $U_N(\alpha)=\{ \{ n\alpha \}, n\in\inter{0}{N} \} $ for any integer $N$. \\
  
%   \esp In the following paragraphs, we will consider, with the order point of view, the Kronecker sequences $( \{ n\alpha \} )_n$ and the more general sequences $(\{ n\alpha + \beta \})_n$, for reals $\alpha$ and $\beta$. But, we begin with some well known results about semi-convergents and best rational in an interval, that will  be useful in the other sections. 

   \subsection{a one-page proof of the "three distance theorem"}

 \esp  In this section, we will be interested in lengths of subdivisions of $[0,1]$ by finite sets $\{ \{ k\alpha\}, k\in \inter{1}{N-1}\}$, where $\alpha$ is a real in $[0,1[$ and $N$ a positive integer. \\
 \esp Let us remark that, if we consider subdivisions of the circle $S^1$, that is to say of $\R/\Z$, then their lengths are invariant by translations. In that case, subdivisions by sets like $\{ \{ k\alpha+\beta\}, k\in \inter{0}{N-1}\}$ are the same, from a metric point of view, for all real $\beta$. \\
 
  \esp The well known 3 distance theorem ( see \textbf{[7]}) claims that these subdivisions are quite simple : they all contains at most 3 different lengths, one being the sum of the others : \\

 \esp Let $\alpha$ be a real in $[0,1[$, with CFE $[a_k]_{k}$. We denote, as usual, $p_n/q_n$ the reduced fraction of the convergent $[a_0,\cdots,a_n]$ and $\delta_n=(-1)^n(\alpha q_n-p_n)$. We remind that $(\delta_n)_n$ is a positive and decreasing sequence that converges towards $0$ ( if $\alpha$ is irrational).\\ 
 \esp  Let $N$ be a positive integer. If $\alpha$ is rational, we suppose that $N\leqslant q$, where $q$ is the denominator of the reduced fraction of $\alpha$. So, the set $\{ \{ k\alpha \},k\in\inter{0}{N-1} \}$ contains exactly $N$ elements.
  
 \begin{montheo}[ 3 distance theorem] .\\
the set $\{ \{ k\alpha \},k\in\inter{1}{N-1} \}$ divides $[0,1]$ into $N$ intervals of length taking at most 3 values, one being the sum of the others. \\
 \esp We can precise a bit : let $s$ be the lowest integer such that $N\leqslant q_s+q_{s-1}$, then :\\
 - if $N=q_s+(1-i)q_{s-1}$, with $i\in\inter{0}{a_s-1}$, the lengths of above intervals take 2 values : 
 $$\delta_s+i\delta_{s-1} \text{ and } \delta_{s-1}$$
 
 - if $N\not =q_s+(1-i)q_{s-1}$,  with $i\in\inter{0}{a_s-1}$, the lengths of above intervals take 3 values : 
 $$\delta_{s-1},\delta_s+i\delta_{s-1} \text{ and }\delta_s+(i+1)\delta_{s-1}$$
 \end{montheo}   
  \begin{mademo}
   \esp  According to propositions 2 and 4,  algorithm 1 and 3, we can write  : $N-1=(n_1,\cdots,n_s)_{\alpha}$, with $n_s\not = 0$. Let denote $(u_j)_{j\in \inter{0}{N-1}}$ the increasing sequence that enumerates our set $\{ \{ k\alpha \},k\in\inter{0}{N-1} \}$. We have $u_0=0$ and denote $u_N=1$. The aim of this result is to prove that $u_j-u_{j-1}$ take at most 3 values, when $j$ ranges over $\inter{1}{N}$.\\
   \esp We will denote $E(N)$ the set of $\alpha$-admissible sequences that are lower or equal, for RLO, than $(n_i)_i$. These sequences are the $\alpha$-numeration of integers of $\inter{0}{N-1}$. Let $k\in \inter{1}{N-1}$, then $k=(k_1,\cdots,k_r)_{\alpha}$, $(k_i)_i \in E(N)$ and $k_r>0$. So, $1\leqslant r \leqslant s$.\\
   \esp We denote $j$ the integer such that $u_j=\{ k\alpha\}$. Then $u_{j-1}=\{ k'\alpha \}$, where $k'=((k'_i)_i)_{\alpha}$ and $(k'_i)_i$ is the predecessor of $(k_i)_i$ in $(E(N),ALO)$. In a similar way $u_{j+1}=\{ k"\alpha \}$, where $k"=((k"_i)_i)_{\alpha}$ and $(k"_i)_i$ is the successor of $(k_i)_i$ in $(E(N),ALO)$.\\
   \esp We will suppose that $s$ is even, because the other case can easily be deduced ( see end of the proof).
   
     $\blacktriangleright$ Case 1 : if $r$ is odd.  
    Then $k'=(k_{[1,r]},1,(\max,0)^{\nu})_{\alpha}$, where $\nu=\frac{s-r-1}{2}$. So, more explicitly :
  %  We define $K'$ :
    $$ (k'_i)_i=(k_1,\cdots,k_{r-1},k_r,1,a_{r+2},0,a_{r+4},0,\cdots,a_{s-1},0)$$
  %  \esp Remark : we have $K'=(k_1,\cdots,k_{r-1},k_r,1)$ if $r=s-1$. Then $K'$ is the predecessor of   $(k_i)_i$ in $(E(N),ALO)$. 
  \esp So :
    $$u_j-u_{j-1}= \delta_{r}- \sum_{i=r+2 ; i \text{ odd}}^{s-1}a_i\delta_{i-1}=\delta_{r}- \sum_{i=r+2 ; i \text{ odd}}^{s-1}(\delta_{i-2}-\delta_i)=\delta_{s-1}$$
  
  $\blacktriangleright$ Case 2 : if $r$ is even. Then, we define $K=(k_{[1,r-1]},k_r+1,(\max,0)^{\nu})$, where  $\nu=\frac{s-r}{2}$.\\
   % $$K=(k_1,\cdots,k_{r-1},k_r+1,a_{r+1},0,a_{r+3},0,\cdots,a_{s-1},0)$$
  $\blacktriangleright\blacktriangleright$ subcase 1 : if $K\in E(N)$, then $K$  is the predecessor of $(k_i)_i$ in $(E(N),ALO)$ and :
  $$u_j-u_{j-1}= \delta_{r-1}- \sum_{i=r+1 ; i \text{ odd}}^{s-1}a_i\delta_{i-1}=\delta_{r-1}- \sum_{i=r+1 ; i \text{ odd}}^{s-1}(\delta_{i-2}-\delta_i)=\delta_{s-1}$$
   $\blacktriangleright\blacktriangleright$ subcase 2 : if $K\not\in E(N)$. We have then 2 subsubcases : \\
    $\blacktriangleright\blacktriangleright\blacktriangleright$ subsubcase 1 : if $r<s$, then $k_r=a_r$. We denote $K'=(k_{[1,r]},0,1,(\max,0)^{\nu})_{\alpha}$, where $\nu=\frac{s-r-2}{2}$.
    %  $$K'=(k_1,\cdots,k_{r},0,1,a_{r+3},0,\cdots,a_{s-1},0)$$
    \esp Then, $K'\in E(N)$ and  $K'$ is the predecessor of $(k_i)_i$ in $(E(N),ALO)$. So :
    $$u_j-u_{j-1}= \delta_{r+1}- \sum_{i=r+3 ; i \text{ odd}}^{s-1}a_i\delta_{i-1}=\delta_{r+1}- \sum_{i=r+3 ; i \text{ odd}}^{s-1}(\delta_{i-2}-\delta_i)=\delta_{s-1}$$
    $\blacktriangleright\blacktriangleright\blacktriangleright$ subsubcase 2 : if $r=s$ then  $k_s=n_s$ or ($k_s=n_s-1$ and $(k_{[1,s-1]})>_R (n_{[1,s-1]})$.\\
    % We have 2 cases : \\
   % $\blacktriangleright\blacktriangleright\blacktriangleright\blacktriangleright$ if $k_{s-1}>0$, then $K"=(k_1,\cdots,k_{s-2},k_{s-1}-1,0)$ is the predecessor of   $(k_i)_i$ in $(E(N),ALO)$ and :
   % $$u_j-u_{j-1}= \delta_{s-2}-k_s\delta_{s-1} $$ 
   %  $\blacktriangleright\blacktriangleright\blacktriangleright\blacktriangleright$ if $k_{s-1}=0$, then
  \esp We denote $t$ the greatest odd integer $i$ such that $k_i>0$. So $k= (k_{[1,t]},(\max,0)^{\nu},k_s)_{\alpha}$, where $\nu=\frac{s-t-1}{2}$. Then, the predecessor of   $(k_i)_i$ in $(E(N),ALO)$ is $(k_{[1,t-1]},k_t-1)$. So :

    $$u_j-u_{j-1}= \delta_{t-1}-k_s\delta_{s-1}- \sum_{i=t+1 ; i \text{ even}}^{s-2}a_i\delta_{i-1}=\delta_{t-1}-k_s\delta_{s-1}- \sum_{i=t+1 ; i \text{ even}}^{s-2}(\delta_{i-2}-\delta_i)=\delta_{s-2}-k_s\delta_{s-1}$$

 N.B : $r=s$ and $k_s=n_s$ is valid for $k=N-1$. But, $r=s,k_s=n_s-1$ and $(k_{[1,s-1]})>_R (n_{[1,s-1]})$ is possible for at least one $k<N$ if and only if $(n_{[1,s-1]})\not = (\max^{s-1})$. That is to say if and only if : $N\not = q_s+q_{s-1}-(a_s-n_s)q_{s-1}=q_{s-2}+(n_s+1)q_{s-1}$.\\
 \esp So, the length $\delta_{s-2}-n_s\delta_{s-1}$ always  occur in our subdivision, but the length $\delta_{s-2}-(n_s-1)\delta_{s-1}$ occur if and only if $N\not =q_{s-2}+(n_s+1)q_{s-1}$. We put $i=a_s-n_s$ and obtain the conditions of Theorem 1.\\    
    
$\blacktriangleright$ Case 3 : the last interval. What about $1-u_j$, where $j=(K)_{\alpha}$ and $K $ is the greatest element of $(E(N),ALO)$ ? Then $K=((\max,0)^{s/2})$, so : 
%$$ K= (a_1,0,a_3,0,\cdots,a_{s-1})$$

$$ 1-u_j= 1- \sum_{i=1 ; i \text{ odd}}^{s-1} a_i\delta_{i-1}=1- \sum_{i=1 ; i \text{ odd}}^{s-1}( \delta_{i-2}-\delta_i)=1-\delta_{-1}+\delta_{s-1} =\delta_{s-1}$$

\esp So, the case $s$ even is proven ! \\

\esp If $s$ is odd, we use similar arguments, replacing "predecessor" by " successor" and "$u_j-u_{j-1}$ " by "$ u_{j+1}-u_j$".  
   \end{mademo}

 \newpage
 
 \subsection{order coincidence of $( \{ n\alpha \} )_n$ and $( \{ n\alpha ' \} )_n$}
   
   $\bullet$ Let $\alpha$ and $\alpha'$ be two different reals in $[0,1)$. We look for the greatest $N$ such that $( \{ n\alpha \})_{n\in \inter{0}{N-1}}$ and $( \{ n\alpha' \})_{n\in \inter{0}{N-1}}$
   are in the same order in the following meaning : 
   $$( \{ n\alpha \})_{n\in I} \text{  is in the same order than } (\{ n\alpha' \} )_{n\in I} \text{ if and only if } ( \forall n,n'\in I, \{ n\alpha \} < \{ n'\alpha \} \Leftrightarrow \{ n\alpha' \} < \{ n' \alpha' \}   )$$ 
   \esp where $I$ is an interval of $\Z$.\\
  \esp This property is related with another one, concerning integral parts :    
    
  \begin{monlem} let $\alpha,\alpha' \in \R$ and $N$ a positive integer. The following assertions are equivalent :  \\
  (i) $( \{ n\alpha \})_{n\in \inter{0}{N-1}} $ and  $ ( \{ n\alpha' \})_{n\in \inter{0}{N-1}} $ are in the same order.\\
  (ii) $\forall n\in \inter{0}{N-1}, \lfloor n\alpha \rfloor = \lfloor n\alpha' \rfloor$
  \end{monlem}
  \begin{mademo} Let $n,n' \in \inter{0}{N-1}$ such that $n<n'$. We denote $d=n'-n\in \inter{0}{N-1}$. Then :
  $$ \lfloor n'\alpha \rfloor=\lfloor d\alpha \rfloor + \lfloor n\alpha \rfloor+\epsilon \esp \text{ where } \epsilon\in \{ 0,1 \}$$
  so :
  $$ \{ n'\alpha \}-\{ n\alpha \}=\{ d\alpha \}-\epsilon$$
  thus, the sign of $ \{ n'\alpha \}-\{ n\alpha \}$ only depends on $\epsilon$. We have the same equalities and remark with $\alpha' $  and $\epsilon'$  instead of $ \alpha$ and $\epsilon$. \\
  $(ii) \Rightarrow (i) : $ suppose that (ii) is true. Then, with above notations, we have $\epsilon=\epsilon'$, so $ \{ n'\alpha \}-\{ n\alpha \}$ and $ \{ n'\alpha' \}-\{ n\alpha' \}$ have the same sign. \\
  $(i) \Rightarrow (ii) : $ suppose that (ii) is false. Then we have an integer $\nu\in \inter{1}{N-1}$ such that :
  $$ \forall k \in \inter{0}{\nu-1}, \lfloor k\alpha \rfloor = \lfloor k\alpha' \rfloor \text{ and } \lfloor \nu\alpha \rfloor\not = \lfloor \nu\alpha' \rfloor$$
  suppose that $\alpha<\alpha'$, then : $ \lfloor \nu\alpha \rfloor < \lfloor \nu\alpha' \rfloor$. If we denote $n'=\nu, n=\nu-1$ and $d=1$, then, with above notations : $\epsilon=0$ and $\epsilon'=1$, so $ \{ n'\alpha \}-\{ n\alpha \}$ and $ \{ n'\alpha' \}-\{ n\alpha' \}$ do not have the same sign.   
    \end{mademo}  
  
  $\bullet$ Suppose that $\alpha$ is a real and $p/q$ is a convergent of $\alpha$. We claim that : \\
% (i) $( \{ n\alpha \})_{n\in \inter{0}{q-1}} $ and  $ ( \{ np/q \})_{n\in \inter{0}{q-1}} $ are in the same order.\\
  $$\forall n\in \inter{0}{q-1}, \esp \lfloor n\alpha \rfloor = \left\lfloor \frac{np}{q} \right\rfloor$$

 \esp Indeed : $\left|\alpha-\frac{p}{q}\right|<\frac{1}{q^2}$, so : $ \forall n\in \inter{1}{q-1}, \left|n\alpha-n\frac{p}{q}\right|<\frac{1}{q}  $. But, $\{ n\frac{p}{q}\}\in [\frac{1}{q},1-\frac{1}{q}]$, since $p$ and $q$ are coprime, so  $\lfloor n\alpha \rfloor = \left\lfloor \frac{np}{q} \right\rfloor$.\\

 $\bullet$ Is this result still valid for semi-convergents instead of convergents ? for other reduced rationals ? The following result gives the answer...and a bit more.

   \begin{maprop}  .\\
   (i) let $\alpha$ and $\alpha'$ be two reals such that $0<\alpha<\alpha'<1$. We denote $\gamma$ the best rational in $]\alpha, \alpha']$ and $q$ the denominator of its reduced fraction. Then    $$ q = \max \{ N\in\N, \forall n\in \inter{0}{N-1}, \lfloor n\alpha \rfloor = \lfloor n\alpha' \rfloor \}$$  
  (ii) let $\alpha$ be a real in $[0,1)$ and $p/q$ a reduced fraction, with $q\in\N^*$, such that $\alpha$ is not the nearest left  strict convergent of $p/q$. 
   $$ p/q \text{ is a semi-convergent of } \alpha \Leftrightarrow  \forall k\in \inter{0}{q-1}, \lfloor k\alpha \rfloor = \lfloor kp/q \rfloor$$ 
   
     \end{maprop}
  
%  \textbf{Remark : } if $\alpha$ is not a semi-convergent of $\alpha'$ ( in particular if $\alpha$ is irrational), then the best rational in $]\alpha, \alpha']$ is the best rational in $[\alpha, \alpha']$.
   
%  \textbf{Remark 1 } " strict" means here that $\alpha\not = p/q$.\\

  \textbf{Remark  : }   for a positive integer $n$, we have $ \lfloor n\alpha \rfloor < \lfloor n\alpha' \rfloor$ if and only if there exists an integer $p$ such that $\alpha < p/n \leqslant \alpha'$.\\  
   
   \begin{mademo}
   (i) is a consequence of the remark.  \\
   (ii) the best rational in $\overset{\longleftrightarrow}{[\alpha,p/q]}$ is the common semi-convergent of $\alpha$ and $p/q$, that has the greatest denominator ( see Proposition 1 (iii)). But, semi-convergents of $p/q$ are either $p/q$ or $p'/q'$ where $p',q'$ are integers such that $1\leqslant q'<q$. So, we have two cases. \\
   \esp If $p/q$ is a semi-convergent of $\alpha$, then there are no integers $a,b$ such that $b\in\inter{1}{q-1}$ and $\alpha < a/b \leqslant p/q$ or $p/q<a/b \leqslant \alpha$. The previous remark implies $\Rightarrow$ of (ii). \\
  \esp If $p/q$ is not a semi-convergent of $\alpha$, then the best rational in $\overset{\longleftrightarrow}{[\alpha,p/q]}$ is $p'/q'$ with $p',q'$ two integers such that $0<q'<q$. If $p/q<\alpha$ then $p/q<p'/q'\leqslant \alpha$ and we use remark 2. Else, since $\alpha$ is  not the nearest left  strict convergent of $p/q$, we have $p",q"$ two integers such that $\alpha<p"/q"<p/q$ and $0<q"<q$.  We conclude with remark 2.
   \end{mademo}

$\bullet$ We also have direct consequences for sums of $\lfloor k\alpha \rfloor $ and $\{ k\alpha \}$ : we will denote
$$ \forall n\in\N, \forall x\in\R, \esp I_n(x)=\sum_{k=0}^{n-1}\lfloor kx \rfloor \esp ; \esp F_n(x)=\sum_{k=0}^{n-1}\{ kx \}$$
\esp Obviously, $F_n$ is 1-periodic, $I_n$ is non decreasing and  : 
$$ \forall n\in\N, \forall x\in\R, \esp I_n(x)+F_n(x)=\frac{n(n-1)x}{2} \esp (1)$$
\esp Moreover, let $p,n$ be 2 positive integers and $d=\gcd(p,n)$. We denote $n'=n/d$ and $p'=p/d$. Then $n'$ and $p'$ are coprime, so $\left\{ \left\{ \frac{kp'}{n'} \right\} ,k\in\inter{0}{n'-1} \right \}= \left \{ \frac{j}{n'} , j\in\inter{0}{n'-1}\right  \} $. So, we have, since $\left( \left\{ \frac{kp'}{n'} \right\}\right)_k$ is $n'$-periodic :
$$ \forall p,n \in\N^*, \esp F_n\left(\frac{p}{n}\right)= \frac{n-\gcd(p,n)}{2} \esp (2)$$
  \esp We also have, for two reals $x$ and $x'$ : 
  $$ I_n(x)=I_n(x') \Leftrightarrow \forall k\in\inter{0}{n-1}, \lfloor kx \rfloor=\lfloor kx' \rfloor$$
  \esp So, Proposition 8 gives : $I_n(x)=I_n(x')$ if and only if $n$ is lower or equal to the denominator of the reduced best rational in $]x,x']$, if $x<x'$. \\
 \esp In  \textbf{[2]}, we can find  an expression of $I_n(x)$ and $F_n(x)$ in terms of the Ostrowski $x$-numeration of $n$. In what follows,  we restrict ourselves to a special case :
  
  \begin{moncoro}
  Let $\alpha$ be a real and $p/q$ a fraction of integers, such that $\alpha$ is not the nearest left  strict convergent of $p/q$. 
    $$ \frac{p}{q} \text{ a  reduced semi-convergent of } \alpha \esp \Leftrightarrow \esp \sum_{k=0}^{q-1} \lfloor k\alpha \rfloor = \frac{(p-1)(q-1)}{2}$$
  
  \end{moncoro} 
   \begin{mademo}
    direct consequence (1),(2) and Proposition 8 (ii). 
      \end{mademo}
      
 % \textbf{Remark : } this result is also true for $\alpha \in\R$, see Annex in 6.2.\\    

  \textbf{Remark : } we deduce an expression of the mean value of $(\{ k \alpha \})_{1\leqslant k<q}$ if $\frac{p}{q} $ is a reduced semi-convergent of $\alpha$ :
  $$ \frac{1}{q-1}   \sum\limits_{k=1}^{q-1}\{ k\alpha \} = \frac{1}{2}+\frac{q\alpha - p }{2}$$
%  \esp This gives a well known fact : this mean value ( that converges towards $1/2$) changes infinitely many times from position beside its limit $1/2$.

%  \vspace*{1cm}
\newpage

   \subsection{best left or right $\alpha$-approximation of a real in $[0,1[$} 
   
   \esp Let $\alpha$ be a real, $[a_k]_{k\in \N^*}$ its CFE and $r=\mu(\alpha)$, the CFE-depth of $\alpha$. So, we denote $[a_0,a_1,\cdots,a_r,1]$ the CFE of $\alpha$ if $\alpha$ is rational. We also denote $(p_n/q_n)_n$ the usual sequence of convergents of $\alpha$. We consider points of $\R^2$ with the product order : $(x,y)\leqslant (x',y')$ if and only if $x\leqslant x'$ and $y\leqslant y'$. \\
  \esp We recall some notations mentioned at  3.2 : 
 for any $x$ in $\R, ||x||$, the distance of $x$ to $\Z$. We also have : $||x||=\min(\{ x \},\{ -x \})$.

  \begin{madef} [ best $\alpha$-approximation of a real].\\
   let $\alpha$ and $\beta$ be two reals in $[0,1[$ and $n$ a non negative integer. \\
  $\triangleright $ \esp $\{ n\alpha \}$ is a \emph{best $\alpha$-approximation} of $\beta$ if and only if :
  $$ \forall k \in\inter{0}{n-1}, \esp ||  n \alpha-\beta||< ||  k \alpha-\beta ||$$
   $\triangleright$   \esp $\{ n\alpha \}$ is a \emph{ best right ( resp. left) $\alpha$-approximation} of $\beta$ if and only if :
  $$ \forall k \in\inter{0}{n-1}, \esp \{   n \alpha-\beta \}<  \{ k \alpha -\beta \} \esp ( \text{ resp. } \{   \beta-n \alpha \}<  \{ \beta-k \alpha  \}) $$
  
  \end{madef}

  \textbf{Remarks : } we could also consider approximations of $\beta$ by $n\alpha$ mod 1, for negative integers $n$. \\
   \esp Best sided $\alpha$-approximations of a real are easier to describe than best $\alpha$-approximations. But, there is a simple relation : a best $\alpha$-approximation is also a best right or left $\alpha$-approximation of $\beta$.\\
 \esp  First, we remark that these notions are closely related to minimal points in $\R^2$  of  sequences $(\{ n\alpha -\beta \},n)_{n\in\N}$ and $(\{ \beta-n\alpha  \},n)_{n\in\N}$ : best right ( resp. left) $\alpha$-approximations of $\beta$ are obtained for the values of $n$ such that $(\{ n\alpha -\beta \},n)$ ( resp. $(\{ \beta-n\alpha  \},n)$) is a minimal point of the sequence $(\{ k\alpha -\beta \},k)_{k\in\N}$ ( resp. $(\{ \beta-k\alpha  \},k)_{k\in\N}$). \\
 \esp Moreover : 
  
  $$ \forall x\in \R, \{ x-\beta \} = \begin{cases} \{ x \} -\beta \in [0,1-\beta [ \text{ if } \{ x \} \geqslant \beta \\
              \{ x \} + 1- \beta \in [ 1- \beta , 1 [ \text{ if } \{ x \} < \beta \end{cases} $$
 
 \esp Finally : $(1-\beta,0)$ is a trivial minimal point of $(\{ n\alpha -\beta \},n)_{n\in\N}$, so the other minimal points must verify $\{ n\alpha \}\geqslant \beta$.

  \begin{maprop} [best right ( positive) $\alpha$-approximations].\\
 $\triangleright$ \textbf{ Case 1} : $\alpha$ is rational and $[0,a_1,\cdots,a_r,1]$ is its CFE. We suppose that $\beta\in \{ \{ n\alpha \}, n\in\N \}$ and denote $(b_1,b_2,\cdots,b_r)$ the $\alpha$-numeration of $\beta$ ( see 2.2).  \\
 \esp Best right ( positive) $\alpha$-approximations of $\beta$ are the $\{ n\alpha \}$  for  $n=0$, for  $n=\sum\limits_{i=1}^rb_iq_{i-1}$ and  for the following $n$  :

  $$ n= \sum_{i=1}^{2k-1}b_iq_{i-1} + j q_{2k-1} \esp ; \esp j\in \inter{0}{b_{2k}-1} \esp ; \esp  k \in\inter{1}{\lfloor r/2 \rfloor}$$

$\triangleright$ \textbf{Case 2 : } if $\alpha$ is irrational and $[a_k]_{k\in\N}$ is its CFE. Let $\beta$ be a real in $[0,1[$ and $(b_k)_{k\in\N^*}$ its $\alpha$-numeration. ( see 2.3)\\
 \esp    Best right ( positive) $\alpha$-approximations of $\beta$ are the $\{ n\alpha \}$  for $n=0$, for $n=\sum\limits_{i=1}^sb_iq_{i-1}$, if $b_k=0$ for all integer $k>s$, and for the following $n$ :

  $$ n= \sum_{i=1}^{2k-1}b_iq_{i-1} + j q_{2k-1} \esp ; \esp j\in \inter{0}{b_{2k}-1} \esp ; \esp  k\in \N^* $$ 
  \end{maprop}

% \textbf{Remark 1 : } these formulae are still true if we replace, in the range of $k$, $t$ by $1$. \\
  
% \textbf{Remark 2 : } if $\beta >0$ and if we want to restrict the range of $n$ in  $\inter{0}{\nu}$, where $\nu\in \inter{0}{q-1}$ : the minimal points in $\R^2$, with the product order, of the sequence $(\{ n\alpha -\beta \},n)_{0\leqslant n \leqslant \nu }$ are those for $n$ ( with $\Psi_{\alpha}$-numeration) : \\
% The same as above that satisfy : $(d_i)_{i\in\inter{1}{r}} \leqslant_R (\nu_i)_{i\in\inter{1}{r}} =\Psi_{\alpha}^{-1}(\nu)$.\\

  \begin{mademo}
  \esp We denote $t=\min(\{ i, b_{2i}\not = 0 \})$, except if all $b_{2i}$ are null : we then denote $t$ the greatest integer $i$ such that $b_{2i-1}\not =0$: so we have, in that case, $b=(\max,0)^t$. Then, for all cases  ( see definition of $E_{(\alpha)}$), we have :
  $$ b_{2t-1}\not = 0 \esp ; \esp b=((\max,0)^{t-1},b_{[2t-1,\infty]})$$
  Following last remarks above Proposition 9, we need the $\alpha$-numeration, say $\nu$, of the least integer $n$ such that $\{ n\alpha \} \geqslant \beta$. According to  Proposition 2,3,4, it is the minimum of elements $d$ of $E_{(\alpha)}$ for RLO, such that $d \geqslant_A b$. We claim that $\nu= b_{[1,2t-1]}$.  Indeed, the condition $d \geqslant_A b$ implies that 
  $$ d_{[1,2t-2]}=(max,0)^{t-1}= b_{[1,2t-2]} \esp \text{ and  } \esp d_{2t-1}\geqslant b_{2t-1}$$
  \esp  But, $b_{[1,2t-1]}$ is minimal ( for RLO) among these one and satisfies $\nu\geqslant_A b$.\\
  \esp Now, if we denote $n_1 =\Psi_{\alpha}(\nu)$ this least integer $n$ such that $\{ n\alpha \} \geqslant \beta$, then : 
  $$ \forall n<n_1, \esp  \{ n\alpha - \beta \}\in [ 1- \beta,1[ \esp ; \esp \{ n_1\alpha - \beta \}\in [ 0,1- \beta[ $$
  \esp So, for the product order in $\Z^2$ :
  $$ \forall n\in \inter{1}{n_1-1}, \esp (1-\beta,0)< (\{ n\alpha - \beta \},n)$$
  \esp Hence, no points $(\{ n\alpha - \beta \},n)$ is minimal, for $n\in \inter{1}{n_1-1}$.\\
  \esp If $b_k=0$ for all integer $k\geqslant 2t$, then  $\nu=b$ and $ \{ n_1\alpha - \beta \}=0$, so this gives the only minimal point ( with $n=0$). \\
  \esp For the other cases :  if $n\geqslant n_1$, let denote $d$ its $\alpha$-numeration. Then, the minimality condition for $(\{ n\alpha - \beta \},n)$ is equivalent to : $d\geqslant_A b$ and $ d$ is minimal among these ( elements of $E_{\alpha}$ greater than $b$ for ALO) for the product of orders (ALO,RLO).  \\
  \esp Of course, $\nu$ is the first ( for RLO) of these minimal ( for (ALO,RLO)) elements. The next one ( for RLO) must satisfy : $d<_A \nu $ and $d$ is minimal for RLO : it gives the successive   $(b_{[1,2t-1]},j) , j\in\inter{0}{b_{2t}-1}$ and then $(b_{[1,2t+1]},j),  j\in\inter{0}{b_{2t+2}-1}$ if $b_{2t+2}\not = 0$ ( but this is still true, if $b_{2t+2}=0$ !), and so on...% the arguments are the same as for the proof of Lemma 15.
  \end{mademo}
  
%\newpage

 $\bullet$  we have a similar result for best left ( positive) $\alpha$-approximations :
 
   \begin{maprop} [best left ( positive) $\alpha$-approximations].\\
$\triangleright$ \textbf{ Case 1} : $\alpha$ is rational and $[0,a_1,\cdots,a_r,1]$ is its CFE. We suppose that $\beta\in \{ \{ n\alpha \}, n\in\N \}$ and denote $(b_1,b_2,\cdots,b_r)$ the $\alpha$-numeration of $\beta$.  \\
 \esp  Best left ( positive) $\alpha$-approximations of $\beta$ are the $\{ n\alpha \}$ for   $n=\sum\limits_{i=1}^rb_iq_{i-1}$ and  for the following $n$  :

  $$ n= \sum_{i=1}^{2k}b_iq_{i-1} + j q_{2k} \esp ; \esp j\in \inter{0}{b_{2k+1}-1} \esp ; \esp  k \in \inter{0}{\lfloor (r-1)/2 \rfloor} $$

 $\triangleright$ \textbf{Case 2 : } $\alpha$ is an irrational and $[a_k]_{k\in\N}$ is its CFE. Let $\beta$ be a real in $[0,1[$ and $(b_k)_{k\in\N^*}$ its $\alpha$-numeration. \\
 \esp    Best left ( positive) $\alpha$-approximations of $\beta$ are the $\{ n\alpha \}$  for  $n=\sum\limits_{i=1}^sb_iq_{i-1}$, if $b_k=0$ for all integer $k>s$, and  the following $n$ :
 $$ n= \sum_{i=1}^{2k}b_iq_{i-1} + j q_{2k} \esp ; \esp j\in \inter{0}{b_{2k+1}-1} \esp ; \esp  k \in \N $$
 
  \end{maprop}
  \begin{mademo} the proof is similar to those of previous Proposition.
  \end{mademo}

  \newpage

   \subsection{measure of repartition of $(\{ k\alpha \})_{0\leqslant k<\nu}$}

   $\bullet$ If  $\alpha$ is an irrational, we know that the sequence of probability measures $(\mu_n)_n$ defined as below converges ( for weak-star topology) to the Lebesgue measure. 
   $$ \forall \nu\in\N^*, \esp \mu_{\nu}=\frac{1}{\nu}\sum_{k=0}^{\nu-1} D_{\{ k\alpha \}}$$
   \esp where $D_x$ is the Dirac-measure in $x$. \\
   \esp Can we precise these measures  ? That is the aim of the following study. It is sufficient to give an expression of $\mu_{\nu}([0,\beta[)$, where $\beta$ is any real of $[0,1[$. So, we want to count integers $k$ in $\inter{0}{\nu-1}$, such that, given a real $\beta$ in $[0,1[$, we have $\{ k\alpha \}<\beta$. \\

   $\bullet$ Another approach of this question is the following : note $L$ the lattice in $\R^2$ generated by $(1,0)$ and $(\alpha,1)$. What is the cardinality of $L\cap R$, if $R$ is the rectangle : $R=[0,\beta[\times [0,\nu[$ ? \\

   $\bullet$ For two reals $\alpha$ and $\beta$ in $[0,1[$ and for a positive integer $\nu$, we denote $n=(n_k)_k$ and $b=(b_k)_k$ the respective $\alpha$-numeration of $\nu$ and $\beta$.  We denote $\sigma$ the usual shift on sequences. We will also use the two total orders on finite sequences of reals : RLO, denoted $\leqslant_R$ and ALO, denoted $\leqslant_A$ ( see 1.2 and 2.3). \\
 \esp We also denote :
   $$ N(\alpha,\beta,\nu)= \{ k\in \inter{0}{\nu-1}, \{ k\alpha \}< \beta \}\esp ; \esp E(\alpha,\beta,\nu)= \{ d\in E_{(\alpha)}, d<_R n \text{ and } d<_A b \} $$ 
   \esp With the results of section 2.3. we can claim that : $\Psi_{\alpha}$ gives a one to one correspondance between $N(\alpha,\beta,\nu)$ and $E(\alpha,\beta,\nu)$. We will denote $ C(\alpha,\beta,\nu)$ the cardinality of these finite sets.\\
   \esp We will denote $\alpha=[a_k]_{k\in\N}$ the CFE of $\alpha$ ( with $a_0=0$) and $r$ the CFE depth of $\alpha$ ( $r=+\infty$ if and only if $\alpha$ is irrational). We suppose $\nu\leqslant q$ if $\alpha$ is a rational and $p/q$ is a reduced fraction that represents $\alpha$. As in section 3.3, we use the following notations : 
    $$ \alpha_0=\alpha \esp ; \esp \forall k\in \inter{1}{r},\alpha_k=\left\{ \frac{1}{\alpha_{k-1}}\right\}$$
 $$ \nu_0=\nu \esp ; \esp \forall k\in \inter{1}{r-2}, \nu_k= \begin{cases} \lfloor \nu_{k-1}\alpha_{k-1}\rfloor  \text{ if } n_{k}\not = 0 \text{ or } n_{k+1}=0 \\ 
   \lfloor \nu_{k-1}\alpha_{k-1}\rfloor +1 \text{ else } \end{cases}  $$
 $$ \beta_0=\beta \esp ; \esp \forall k\in \inter{1}{r}, \beta_k= \frac{1}{\alpha_{k-1}}(b_k\alpha_{k-1}-\beta_{k-1})$$

  \textbf{Remark 1 : } 
  $$ d\in E_{(\alpha)} \Leftrightarrow d=(0) \text{ or } \begin{cases} d_1\in \inter{1}{a_1} \\ \sigma(d)\in E_{(\alpha_1)} \end{cases} \text{ or }  \begin{cases} d_1=a_1  \\ d_2=0  \\ \sigma^2(d)\in E_{(\alpha_2)}\backslash \{ (0) \}  \end{cases}$$
  \esp These three cases are exclusive.\\
  % and the third one is not possible if $b_2\not = 0$, if we suppose that $d<_A b$. \\
  
  \textbf{Remark 2 : }  let $ d\in E_{(\alpha)}$, then :
   $$ d<_R   n \Leftrightarrow  \sigma(d)<_R \sigma(n)  \text{ or } ( \sigma(d)=\sigma(n) \text{ and } d_1<n_1 )$$
   $$ d <_A b \Leftrightarrow d_1<b_1 \text{ or } ( d_1=b_1 \text{ and } \sigma(b)<_A \sigma(d) )$$

 \vspace*{0.5cm}

   \begin{maprop} we denote $n=(n_k)_k$ the $\alpha$-numeration of $\nu$ and $b=(b_k)_k$ the $\alpha$-numeration of $\beta$. We denote $s$ the minimum of the lengths of $n$ and $b$, when we drop the eventual infinite " $0$-tail". So, $n_s$ or $b_s$ is not null, but $\sigma^s(n)$ or $\sigma^s(b)$ is the null sequence.\\
   $$ C(\alpha,\beta,\nu)= \sum_{i=1}^{s } (-1)^{i-1}[  b_i \nu_i +\tau_i+\epsilon_i-\epsilon'_i]$$

$\tau_i= \begin{cases} 1 \text{ if } n_in_{i+1}=0 \text{ and }  \sigma^i(n)\not = (0) \\ \min(b_i,n_i) \text{ else} \end{cases} $ \\
 $\epsilon_i= \begin{cases} 1 \text{ if }  b_i< n_i \text{ and } \sigma^i(b) <_A \sigma^i(n) \\ 0 \text{ else } \end{cases}$ 
$\epsilon'_i= \begin{cases} 1 \text{ if }   \sigma^i(b) <_R \sigma^i(n) \\ 0 \text{ else } \end{cases}$
     
   \end{maprop}

   \begin{mademo} we want to enumerate sequences $d$ of $E_{(\alpha)}$ such that $d<_R n$ and $d<_A b$. We will consider several cases and subcases, depending on the cancellation of the $b_i$ and $n_i$...\\
   \esp First, we remark that $b_1>0$ and $n_1>0$, for we can suppose that $b\not = (0)$ and $n\not = (0)$. \\
    $\blacktriangleright $ Case 1 : $b_2>0$.\\
      $\blacktriangleright \blacktriangleright $ subcase 1 : $n_2>0$ or $\sigma(n)=(0)$. Let us count sequences $d$ as follows : \\
      --- if $d_1=0$, then $d=(0)\in E(\alpha,\beta,\nu)$ for $n\not = (0)$ and $b\not = (0)$ : 1 sequence. \\
      ---  if $0<d_1<b_1$. Then $d<_A b$. So $d\in E(\alpha,\beta,\nu)$ if and only if $d<_R n$. \\
       ------ if $ \sigma(d)=\sigma(n) \text{ and } d_1<n_1$, this gives, exactly $\min(b_1,n_1)-1$ sequences $d$. \\
    ------ if $  \sigma(d)<_R \sigma(n) $, this gives, for every $d_1\in \inter{1}{b_1-1}$,   $\nu_1$ possible sequences $d$, according to Lemma 6 : so, we have $\nu_1(b_1-1)$ sequences $d$ for this subcase. \\
  --- if  $ d_1=b_1$, then, $d\in  E(\alpha,\beta,\nu)$ if and only if $d<_R n$ and $\sigma(d)>_A \sigma(b)$. \\
  ------ if $\sigma(d)=\sigma(n)$, this gives a unique sequence $d=(b_1,\sigma(n))$ if and only if $b_1<n_1$ and $\sigma(n)>_A\sigma(b)$ ( because $n_2\not = 0$, so $(b_1)\sqcup \sigma(n) \in E_{(\alpha)}$)  and no sequences $d$ else. This gives $\epsilon_1$ sequences.\\
   ------ if $\sigma(d)<_R \sigma(n)$. Since $d_1$ is fixed ( $d_1=b_1)$, counting these sequences is the same, according to Lemma 6, as counting sequences $u$ of $E_{(\alpha_1)}$ ( since $d_2$ can not be null if $\sigma(d) >_A \sigma(b)$) such that $u <_R \sigma(n)$ and $u >_A \sigma(b)$.  But,$\sigma(n)$ is the $\alpha_1$- numeration of $\nu_1$ ( see Lemma 7) and $\sigma(b)$ is the $\alpha_1$-numeration of $\beta_1$ ( see Lemma 8).  So,  we obtain $\nu_1-C(\alpha_1,\beta_1,\nu_1)-\epsilon'_1$ sequences $d$ for this subcase, where $\epsilon'_1=1$ if  and only if $u$ can be equal to $\sigma(b)$, so if and only if  $\sigma(b)<_R \sigma(n)$ and $0$ else. \\  
    
   \esp If we summarize this subcase, we obtain :
   $$ C(\alpha_0,\beta_0,\nu_0)=\nu_1b_1+\min(b_1,n_1)+\epsilon_1-\epsilon'_1-C(\alpha_1,\beta_1,\nu_1)$$

   $\blacktriangleright\blacktriangleright$ subcase 2 : $n_2=0$ and $\sigma(n)\not = (0)$.\\
   --- if $d_1=0$ : 1 sequence for $d=(0)$.\\
   --- if $0<d_1<b_1$, this is the same count as in the previous subcase, except that :  we have $\sigma(n)=(0,n_{[3,\infty]})$ with $n_3\not = 0$, so $\sigma(n)$ is not a possible value for $\sigma(d)\in E_{(\alpha_1)}$ if $d_1<b_1 $ ( for $b_1\leqslant a_1$). So, we must replace $\min(b_1,n_1)$ by $1$  : this is the role of $\tau_1$.  Furthermore, the condition $  \sigma(d)<_R \sigma(n) $ is equivalent to $\sigma(d)<_R (1,n_{[3,\infty]})=(1, \sigma^2(n))$ that is  the $\alpha_1$-numeration of $\nu_1$ : so this gives $\tau_1-1+\nu_1(b_1-1)$ sequences.\\
   --- if $d_1=b_1$, we have $n_2=0$, so $u<_R\sigma(n)$ is equivalent to $ u<_R (1,\sigma^2(n))$ and $(1,\sigma^2(n))$ is the $\alpha_1$- numeration of $\nu_1$. As above,  we obtain $\nu_1-C(\alpha_1,\beta_1,\nu_1)-\epsilon'_1$ sequences $d$ for this subcase. Now, with all previous arguments, we obtain $C(\alpha_1,\beta_1,\nu_1)=\nu_2b_2+\tau_2+\epsilon_2-\epsilon'_2-C(\alpha_2,\beta_2,\nu_2)$, but $\nu_1=(1,\sigma^2(n))_{\alpha_1}$ and $n_2=0$, so we must replace $0$ by $1$ for the value of $n_2$ in the formula for $\tau_2$ and $\epsilon_2$. But, it does not change the result, for $b_2\not = 0$ ! At the end, $\nu_2=(\sigma^2(n))_{\alpha_2}$, so the induction goes on. \\
   
    \esp If we summarize this subcase, we obtain ( here $\epsilon_1=0$) : 
   $$ C(\alpha_0,\beta_0,\nu_0)=\nu_1b_1+\tau_1+\epsilon_1-\epsilon'_1-(\nu_2b_2+\tau_2+\epsilon_2-\epsilon'_2)+ C(\alpha_2,\beta_2,\nu_2)$$

  $\blacktriangleright$ Case 2 : if $b_2=0$ and $\sigma(b)\not = 0$. Then $b_1=a_1$.  We can copy all arguments given in Case 1, except if $d_1=b_1$ and $\sigma(d)<_R\sigma(n)$ : indeed, $\sigma(b)$ is not the $\alpha_1$-numeration of $\beta_1$ ( for $\beta_1<0$ and $\sigma(b)\not \in E_{\alpha_1}$). But, $\sigma^2(b)$ is the $\alpha_2$-numeration of $\beta_2$ ( see Lemma ...). So, we must look for a formula between $C(\alpha_0,\beta_0,\nu_0)$ and $C(\alpha_2,\beta_2,\nu_2)$.  Moreover,  $\sigma(d) >_A \sigma(b)$ if and only if $d_2>0$ or ($d_2=0$ and $\sigma^2(d)<_A \sigma^2(b)$). \\
 ---  if $d_2>0$, then counting these sequences is the same as counting sequences $d$ such that $d_2>0$ and $\sigma(d)<_R\sigma(n)$, so counting sequences $u\in E_{(\alpha_1)}$ such that $u\not = (0)$ and $u<_R \sigma(n)$. With the same arguments as in Case 1 ( separating 2 cases : if $n_2$ is null or not), we obtain $\nu_1-1$ such sequences.\\
 --- if $d_2=0$. We will study 3 subcases, depending on $n_2$ and $n_3$ :\\
    $\blacktriangleright\blacktriangleright$ subcase 1 : if $n_2=0$, then we count sequences $d$ such that $\sigma^2(d)<_R\sigma^2(n)$ and $\sigma^2(d)<_ A\sigma^2(b)$. So, we obtain $C(\alpha_2,\beta_2,\nu_2)$ such sequences, because  $\sigma^2(n)$ and $\sigma^2(b)$ are the $\alpha_2$-numeration of $\nu_2$ and $\beta_2$ respectively. \\
     $\blacktriangleright\blacktriangleright$ subcase 2 : if $n_2\not = 0$ and ( $n_3\not = 0$ or $\sigma^2(n)=(0)$), then we count sequences $d$ such that $\sigma^2(d)\leqslant_R\sigma^2(n)$ :   we obtain $C(\alpha_2,\beta_2,\nu_2)+\epsilon"_1$ such sequences, with $\epsilon"_1=1$ if $\sigma^2(n)<_A\sigma^2(b)$, $\epsilon"_1=0$ else...( $\sigma^2(n)$ and $\sigma^2(b)$ are still the $\alpha_2$-numeration of $\nu_2$ and $\beta_2$ respectively). \\
    $\blacktriangleright\blacktriangleright$ subcase 3 : if $n_2\not =0,n_3=0$  and $\sigma^3(n)\not = (0)$, then $\sigma^2(n)$ is not  the $\alpha_2$-numeration of $\nu_2$ : it is $ (1,\sigma^3(n))$. Now, $\sigma^2(d)\leqslant_R\sigma^2(n)$ is equivalent to $\sigma^2(d)<_R (1,\sigma^3(n))$, so we obtain $C(\alpha_2,\beta_2,\nu_2)$ sequences $d$ ( see Lemma 7 again). \\

     \esp If we summarize this case 2 :

$$ C(\alpha_0,\beta_0,\nu_0)=\nu_1b_1+\tau_1+\epsilon_1-1+\epsilon"_1+C(\alpha_2,\beta_2,\nu_2)$$
\esp where $\epsilon"_1=1$ if $n_2\not =0,(n_3\not = 0$ or $\sigma^2(n)=(0)$) and $\sigma^2(n)<_A\sigma^2(b)$. $\epsilon"_1=0$ else. \\

  \esp Now, let us summarize and generalize all cases :   
     
     for all $i\in \inter{1}{s-1}$ :  ( we have $ \sigma^i(b)\not =(0)$)\\
     - if $b_{i+1}>0$, then : $C(\alpha_{i-1},\beta_{i-1},\nu_{i-1})=\nu_ib_i+\tau_i+\epsilon_i-\epsilon'_i-C(\alpha_i,\beta_i,\nu_i)$.\\
     - if $b_{i+1}=0$ and $\sigma^i(b)\not = (0)$, then : $C(\alpha_{i-1},\beta_{i-1},\nu_{i-1})=\nu_ib_i+\tau_i+\epsilon_i+\epsilon"_i+C(\alpha_{i+1},\beta_{i+1},\nu_{i+1})-1$, \\
    \esp  where $\epsilon"_i=1$ if $n_{i+1}\not = 0,(n_{i+2}\not = 0$ or $ \sigma^{i+1}(n)=(0)$) and $\sigma^{i+1}(n)<_A\sigma^{i+1}(b)$ and $0$ else.\\
    \esp  We claim that   :
      $$\epsilon"_i-1=-\epsilon'_i+\epsilon'_{i+1}-\epsilon_{i+1}-\tau_{i+1}$$
     
     --- if $n_{i+1}=0$, then $\epsilon_{i+1}=0,\tau_{i+1}=1$ and $\sigma^i(b)<_R\sigma^i(n) \Leftrightarrow \sigma^{i+1}(b)<_R \sigma^{i+1}(n)$, so $\epsilon'_i=\epsilon'_{i+1}$. Moreover, $\epsilon"_i=0$, so the equality is true. \\
      --- if $n_{i+1}>0$ , then  $\sigma^i(b)<_R\sigma^i(n) \Leftrightarrow \sigma^{i+1}(b)\leqslant_R \sigma^{i+1}(n)$.\\
      \esp If $\sigma^{i+1}(b) = \sigma^{i+1}(n)$, then $\epsilon_{i+1}=0,\tau_{i+1}=0,\epsilon'_i=1, \epsilon'_{i+1}=0$ and $\epsilon"_i=0$,  so the equality is true. \\
      \esp If $\sigma^{i+1}(b)\not = \sigma^{i+1}(n)$, then $\epsilon'_i=\epsilon'_{i+1}$. If $n_{i+2}\not = 0$ or $ \sigma^{i+1}(n)=(0)$ then $\epsilon"_i=1-\epsilon_{i+1}$ and $\tau_{i+1}=0$. Else, $\epsilon_{i+1}=0$ ( for $b_{i+2}\not = 0$),$\epsilon"_i=0$ and $\tau_{i+1}=1$. In both cases, the equality is true.\\

      \esp From this equality, we deduce that : if $b_{i+1}=0$,  and $ \sigma^i(b)\not =(0)$, then
      $$   C(\alpha_{i-1},\beta_{i-1},\nu_{i-1})=\nu_ib_i+\tau_i+\epsilon_i-\epsilon'_i-(\nu_{i+1}b_{i+1}+\tau_{i+1}+\epsilon_{i+1}-\epsilon'_{i+1})+C(\alpha_{i+1},\beta_{i+1},\nu_{i+1})$$
      So, the induction formula for $b_{i+1}>0$ can be generalized to all cases and we conclude with : if $s=1$, then 
    $n=(n_1)$ or $b=(b_1)$. In the first case, $\nu_1=0$ and   $C(\alpha,\beta,\nu)$ counts the $d=(d_1)$ such that $0\leqslant d_1<n_1$ and $d_1<b_1$. So  $C(\alpha,\beta,\nu)=\min(b_1,n_1)=\tau_1$ and $\epsilon_1=0=\epsilon'_1$, since $\sigma(n)=(0)$. In the second case, we have $\sigma(b)=(0)\not = \sigma(n)$. Our former arguments give : $C(\alpha,\beta,\nu)=b_1\nu_1+\tau_1+\epsilon_1-1$ and $\epsilon'_1=1$. This is the initialization of our induction.     
  \end{mademo}

 $\bullet$ We can deduce similar results for conditions with large inequalities instead of strict ones.\\
 \esp For example : if we denote $C'(\alpha,\beta,\nu)=\# \{ k\in\inter{0}{\nu}, \{ k\alpha \}\leqslant \beta \}$, then :
 $$ C'(\alpha,\beta,\nu)=C(\alpha,\beta,\nu)+D$$
\esp where :  
$$D = \textbf{1}_{ n\leqslant_A b}+ \textbf{1}_{ b\leqslant_R n}-\textbf{1}_{ n=b} $$
    \esp Indeed, if we denote $E'(\alpha,\beta,\nu)=\{ d\in E_{(\alpha)}, d\leqslant_R n, d\leqslant_A b \}$, then $C'(\alpha,\beta,\nu)$ is the number of elements of $E'(\alpha,\beta,\nu)$. This set is $E(\alpha,\beta,\nu)$ plus the element $n$ if and only if $n\leqslant_Ab$, plus the element $b$ if and only if $b\leqslant_R n$... if $n=b$, we have to count once this element.
    
  %  \newpage

\newpage

    \section{References}
    
\label{ref:BER}
[1] V Berth\'e : " autour du syst\`eme de num\' eration d'Ostrowski",
Bull. Belg. Math. Soc. 8 (2001), 209-238 \\

\label{ref:BRO}
[2] T.C. Brown and P.J.-S. Shiue : " sums of fractional parts of integer multiples of an irrational", J.Number Theory 50 (1995), 181-192.\\

\label{ref:CAB}
[3] E.Cabanillas : " quotients of numerical semigroups generated by two numbers", hal-02097473 and Arxiv 1904.08240  ( 2019) \\

\label{ref:CAS}
[4] J. W. S. Cassels : " an introduction to Diophantine approximation", Cambridge, Cambridge University Press, 1957.\\

\label{ref:ITO}
[5] S. Ito : "some skew product transformations associated with continued fractions and their invariant measures" , Tokyo J. Math. 9 (1986), 115-133.\\

\label{ref:OST}
[6] Ostrowski : " bemerkungen zur Theorie der Diophantischen Approximationnen I,II" , Abh. Math. Sem Hamburg I ( 1922), 77-98 and 250-251\\

\label{ref:SOS}
[7] V. T. S S\`os : "on the distribution mod 1 of the sequence n $\alpha$", Ann. Univ. Sci. Budapest, Eotvos Sect. Math. 1 (1958), 127-134.

 \end{document}